%______________________________________________________________
%
%   Sharp L1 Stability Estimates 
%               for  
%  Hyperbolic Conservation Laws
%               by 
%    P. Goatin and P.G. LeFloch 
%   To appear in Portugaliae Mathematica (2001)  
%    Version dated September 21, 2000  
%
%______________________________________________________________
\input amstex 
\documentstyle{amsppt}
\magnification 1100
\TagsOnRight
%\NoBlackBoxes 
\hsize=6.6 truein 
%\vsize=8.5 truein 
%______________________________________
 
\define \del {\partial}
\define \RR  {{R\!\!\!\!\!I~}}
\define \RN  {\RR^N}
\define \ve {\varepsilon}
\define \eps {\epsilon}
\define \JJ  {{\Cal J}}
\define \II  {{\Cal I}}
\define \LL  {{\Cal L}}
\define \SS  {{\Cal S}}
\define \FF  {{\Cal F}}
\define \RS  {{\Cal R}}
\define \CC  {{\Cal C}}
\define \ZZ  {{\Cal Z}}
\define \PP  {{\bold P}}
\define \lam {\lambda}
\define \kap {\kappa} 
\define \sgn {\text{sgn }}  
%_____________________________________
\define \refBressan      {1}
\define \refBCP          {2}
\define \refBL           {3}
\define \refBLY          {4} 
\define \refCLone        {5}
\define \refCLtwo        {6}
\define \refDafermosone  {7}
\define \refDafermostwo  {8}
\define \refDafermosthree{9}
\define \refDLM          {10}
\define \refEG           {11} 
\define \refFilippov     {12} 
\define \refGoatinL      {13}
\define \refHL           {14}
\define \refLax          {15} 
\define \refLeFlochone   {16} 
\define \refLeFlochtwo   {17} 
\define \refLeFlochthree {18} 
\define \refLeFlochfour  {19} 
\define \refLL           {20}  
\define \refLX           {21}
\define \refLYone        {22}
\define \refLYtwo        {23}      
\define \refLYthree      {24}  
\define \refVolpert      {25}
%______________________________________________________________
\topmatter
\title 
Sharp~${\bold L^1}$~Stability~Estimates
for
Hyperbolic~Conservation~Laws   
\endtitle 
\author 
Paola Goatin$^{1,2}$ and Philippe G. LeFloch$^1$ 
\endauthor
\rightheadtext{Sharp $L^1$ stability estimates}
\leftheadtext{P. Goatin and P.G. LeFloch}
\footnote""{ \, Published in: Port. Math. 58 (2001), 1--44.
\newline  
\indent \, $^1$ \,
% Centre de Math\'ematiques Appliqu\'ees \&
%Centre National de la Recherche Scientifique, U.M.R.~7641,
% Former address: Ecole Polytechnique, France.
%91128 Palaiseau Cedex, France.
%E-mail address: lefloch\@cmap.polytechnique.fr. 
Current address:  Laboratoire Jacques-Lions Lions \& Centre National de la Recherche Scientifique (CNRS), 
Ê      Universit\'e Pierre et Marie Curie (Paris 6), 
 Ê     4 Place Jussieu, 75252 Paris, France. 
 {\it E-mail:} pgLeFloch\@gmail.com. 
\newline 
\indent \, $^2$ \, S.I.S.S.A., Via Beirut 4, Trieste 34014, Italy. 
%E-mail:goatin\@sissa.it, goatin\@cmap.polytechnique.fr. 
\hfill \hskip6.cm}  
\abstract  
We introduce a generalization of Liu-Yang's weighted norm 
to linear and to nonlinear hyperbolic equations.
Following an approach due to the second author for piecewise constant solutions, 
we establish sharp $L^1$ continuous dependence estimates for general 
solutions of bounded variation. Two different strategies are successfully 
investigated. On one hand, we justify passing to the limit in an 
$L^1$ estimate valid for piecewise constant wave-front tracking approximations. 
On the other hand, we use the technique of 
generalized characteristics and, following closely an approach by Dafermos, 
we derive the sharp $L^1$ estimate directly from the equation.  
\endabstract
\endtopmatter
\document
%=============================================================== 
\heading{1. Introduction} 
\endheading

We are interested in the continuous dependence of entropy solutions to  
hyperbolic conservation laws 
$$
\del_t u + \del_x f(u) = 0, \qquad u(x,t) \in \RR, \, x \in \RR, \, t>0, 
\tag 1.1 
$$ 
where the flux $f: \RR \to \RR$ is a smooth and convex function. 
After works by Liu and Yang \cite{\refLYone} and Dafermos \cite{\refDafermosthree}, 
we aim at deriving sharp $L^1$ estimates of the form 
$$ 
\|u^{II}(t) - u^I(t)\|_{w(t)} 
+ \int_s^t M(\tau; u^I, u^{II}) \, d \tau 
\leq  
\|u^{II}(s) - u^I(s)\|_{w(s)}, \qquad 0 \leq s \leq t, 
\tag 1.2 
$$ 
for any two entropy solutions of bounded variation $u^I$ and $u^{II}$ of (1.1), 
where $\| . \|_{w(t)}$ is a weighted $L^1$-norm equivalent 
to the standard $L^1$ norm on the real line. In (1.2), the positive 
term $M(\tau; u^I, u^{II})$ is intended to provide a sharp bound on 
the strict decrease of the $L^1$ norm. The estimate with $w \equiv 1$ and $M \equiv 0$ 
is of course well-known.  

Recall that the fundamental issue of the uniqueness and continuous dependence 
for hyperbolic systems of conservation laws was initiated by 
Bressan and his collaborators (see \cite{\refBressan, 
\refBCP} and the references therein). A major contribution came from 
Liu and Yang \cite{\refLYone, \refLYtwo} who introduced a decreasing $L^1$ functional
ensuring (1.2) for scalar conservation laws and systems of two equations. 
This research culminated in papers published simultaneously by 
Bressan, Liu and Yang \cite{\refBLY}, Hu and LeFloch \cite{\refHL}, and 
Liu and Yang \cite{\refLYthree}, which contain particularly simple 
proofs of the continuous dependence of entropy solutions for systems.  

%_________________________________________________________________ 

In the present paper, we restrict attention to scalar conservation laws and, 
following the approach developed by the second author 
(Hu and LeFloch \cite{\refHL} and LeFloch \cite{\refLeFlochthree}), 
we investigate the stability issue from the standpoint of Holmgren's 
and Haar's methods (\cite{\refLX} and the references therein).  
The problem under consideration is (essentially) equivalent 
to showing the uniqueness and $L^1$ stability for the following 
hyperbolic equation with discontinuous coefficient:  
$$ 
\del_t \psi + \del_x \bigl(a\, \psi \bigr) = 0, \qquad \psi(x,t) \in \RR,
\, x \in \RR, \, t>0. 
\tag 1.3
$$
That is, for solutions with bounded variation we aim at deriving
an estimate like 
$$ 
\|\psi(t)\|_{w(t)} + \int_s^t \tilde M(\tau; a, \psi) \, d \tau  
\leq  
\|\psi(s)\|_{w(s)}, \qquad 0 \leq s \leq t. 
\tag 1.4 
$$ 
For the application to (1.1) one should define $a$ by 
$$
a = a(u^I, u^{II}) = {f(u^{II}) - f(u^I) \over u^{II} - u^I}.
\tag 1.5  
$$ 
One may also consider the equation (1.3) for more general coefficients $a$. 

Recall that the existence and uniqueness of solutions 
to the Cauchy problem associated with (1.3) 
was established in LeFloch \cite{\refLeFlochone} in the class of bounded measures, 
under the assumption $a_x \leq E$ for some constant $E$. The latter
holds when $a$ is given by (1.5) (at least when $u^I$ and $u^{II}$ contain 
no rarefaction center on the line $t=0$ which holds ``generically''). 
See also Crasta and LeFloch \cite{\refCLone, \refCLtwo} for further existence results.

It must be observed that 
we restrict attention here to more regular solutions, having bounded 
total variation, as this is natural in view of the application to the conservation 
law (1.1). In this direction, 
recall that an $L^1$ stability result like (1.4) was established in 
\cite{\refHL} (see, therein, Section 5, and our Theorem 2.2 below) 
in the class of piecewise Lipschitz continuous solutions, with $\tilde M \equiv 0$ however.  
This uniqueness and stability result was achieved under the assumption 
that the coefficient $a$ does not contain any rarefaction shock 
(see Section 2 below for the definition). In \cite{\refHL},  
the following essential observation was made:
$$ 
\aligned 
\text{The linearized equation $(1.3)$-$(1.5)$ based on  
two entropy }\\ 
\text{
solutions of $(1.1)$ does not exhibit rarefaction shocks.} 
\endaligned  
\tag 1.6 
$$ 
(This is also true for systems of conservation laws, as far as  
solutions with small amplitude are concerned.)  
One of our aims here is to extend the $L^1$ stability result for (1.3) 
in \cite{\refHL} to arbitrary solutions of bounded variation. 

%_______________________________________________________________ 

The present paper relies also heavily on the contribution by Liu and Yang 
\cite{\refLYone} who, for approximate solutions constructed by the Glimm scheme, 
discovered a weighted norm having a sharp decay of the form 
(1.2). Subsequently, the Liu-Yang's functional was extended by Dafermos 
(\cite{\refDafermosthree}, Chapter 11) 
to arbitrary functions of bounded variation (BV) and, using the notion of 
generalized characteristics,  
Dafermos derived precisely an estimate of the form (1.2) valid for BV solutions. 
 
The aim of this paper is to provide a new derivation and some generalization 
of this $L^1$ functional. 
Toward the derivation of bounds like (1.2) or (1.4), we make the 
following preliminary observations: 
\roster
\item The geometrical properties of the propagating discontinuities in $a$ 
(Lax, fast or slow undercompressive, rarefaction shock, according to the 
terminology in \cite{\refHL}) play an essential role. 
It turns out that the (jump of the) weight $w(x,t)$ should be assigned  
precisely on each {\it undercompressive discontinuity.\/} On the other hand,
Lax discontinuities are very stable and do not require  weight, 
while (in exact entropy soutions) 
rarefaction shocks do not arise, according to (1.6). 
\item Certain (invariance) properties on the coefficient $a$ are necessary 
to define the weight globally in space; see (2.9)-(2.10) in Section 2. 
\item The weight however is far from being unique and we believe that 
this flexibility in choosing the weight may be helpful in certain applications. 
\endroster 

The content of this paper is as follows. 

In Section 2, we consider piecewise constant solutions of (1.3) and introduce 
a class of weighted norms satisfying a sharp bound of the form (1.4). 
See Theorem 2.3 below. All undercompressive and Lax discontinuities contribute 
to the decrease of the $L^1$ norm. For the sake of comparison, we also consider 
the $L^1$ norm without weight; see Theorem 2.2. 

In Section 3, we point out that the setting of Section 2 covers the case of
the conservation law (1.1). 
Passing to the limit in wave front tracking approximations, in Theorem 3.5 
we arrive to the sharp bound (1.2) for general BV solutions. The proof is based 
on fine convergence properties established earlier by Bressan and LeFloch \cite{\refBL} 
and on a technique of stability of nonconservative products developed by 
DalMaso, LeFloch, and Murat \cite{\refDLM} and LeFloch and Liu \cite{\refLL}. 

Next, in Sections 4  and 5 we return to the equation (1.3) studied in Section 2 
but, now, we deal with general BV solutions. We follow closely ideas
developed by Dafermos \cite{\refDafermostwo, \refDafermosthree} for solutions of (1.1), 
and extend them to the linear equation (1.3). 
Using generalized characteristics we establish first a maximum principle in Theorem 4.5.  
Finally, in Theorem 5.1 using the technique of generalized characteristics, 
we establish the sharp $L^1$ stability property (1.4) 
directly, for general BV solutions of (1.3). The result applies in particular to the 
conservation law (1.1) and allows us to recover (1.2).

Throughout the paper, we always assume that all functions of bounded 
variation under consideration are normalized to be defined everywhere 
as right-continuous function. 

%=============================================================== 
\heading{2. Decreasing Weighted Norms for Piecewise Constant Solutions} 
\endheading

Given a piecewise constant function $a: \RR\times \RR_+ \to \RR$, 
let us consider the linear hyperbolic equation 
$$ 
\del_t \psi + \del_x \bigl(a\, \psi \bigr) = 0, \qquad \psi(x,t) \in \RR, 
\tag 2.1 
$$ 
and restrict attention to piecewise constant solutions. 
By definition, the function 
$a$ admits a set of jump points $\JJ(a)$, consisting of finitely 
many straightlines defined on open time intervals, together with 
a finite set of interaction points $\II(a)$,
consisting of the end points of the lines in $\JJ(a)$. The function   
$a$ is constant in each connected component of the complement $\CC(a)$ 
of $\II(a) \cup \JJ(a)$.
At a point $(x,t) \in \JJ(a)$ we denote by $\lam^a=\lam^a(x,t)$ the speed of 
the discontinuity and $a_\pm= a_\pm(x,t)= a(x\pm, t)$ 
the left- and right-hand traces.  
It is tacitly assumed that the discontinuity speeds $\lam^a$ remain uniformly bounded.  
Finally the function is normalized to be right-continuous. 
A similar notation is used for the function $\psi$. 

The geometrical property of the coefficient $a$ play a central role for the 
analysis of (2.1), 
so we recall the following terminology \cite{\refHL}:

\proclaim{Definition 2.1}  
A point $(x,t) \in \JJ(a)$ is called a Lax discontinuity iff 
$$ 
a_-(x,t) > \lam^a(x,t) > a_+(x,t),  
$$ 
a slow undercompressive discontinuity iff 
$$
\lam^a(x,t) \leq \min\bigl(a_-(x,t), a_+(x,t)\bigr), 
$$ 
a fast undercompressive discontinuity iff 
$$
\lam^a(x,t) \geq \max\bigl(a_-(x,t), a_+(x,t)\bigr), 
$$ 
and a rarefaction-shock discontinuity iff 
$$ 
a_-(x,t) < \lam^a(x,t) < a_+(x,t). 
$$ 
\endproclaim 

For each $t>0$, we denote by $\LL(a), \SS(a), \FF(a)$, and $\RS(a)$ 
the set of points $(x,t) \in \JJ(a)$ corresponding to Lax, slow undercompressive, 
fast undercompressive, and rarefaction-shock discontinuities, respectively. 

%_____________________________________________________________________________ 

\proclaim{Theorem 2.2} Consider a piecewise constant speed $a=a(x,t)$. 
Let $\psi$ be any piecewise constant solution of $(2.1)$.  
Then we have for all $0 \leq s \leq t$  
$$
\aligned 
& \|\psi(t)\|_{L^1(\RR)} 
+ 
\int_s^t \sum_{(x,\tau) \in \LL(a)} 
2 \, \bigl(a_-(x,\tau) - \lam^a(x,\tau) \bigr) \, |\psi_-(x,\tau)| \, d\tau \\ 
& = 
\|\psi(s)\|_{L^1(\RR)} + 
\int_s^t \sum_{(x,\tau) \in \RS(a)} 
2 \, \bigl(\lam^a(x,\tau) - a_-(x,\tau) \bigr) \, |\psi_-(x,\tau)| \, d\tau. 
\endaligned 
\tag 2.2  
$$
\endproclaim

In (2.2), the left-hand traces are chosen for definiteness only. 
Indeed it will be noticed in the proof below that 
for all $(x,\tau) \in \LL(a) \cup \RS(a)$ 
$$
\bigl(\lam^a(x,\tau) - a_-(x,\tau) \bigr) \, |\psi_-(x,\tau)| 
= 
- \bigl(\lam^a(x,\tau) - a_+(x,\tau) \bigr) \, |\psi_+(x,\tau)| 
$$
Observe that the Lax discontinuities contribute to the decrease of the $L^1$ norm, 
while the rarefaction-shocks increase it. 
On the other hand, the undercompressive discontinuities don't modify the $L^1$ norm. 
When $a$ contains no rarefaction shocks (this is the case when 
(2.1) is a linearized equation derived from entropy solutions of a conservation 
law, as discovered in Hu and LeFloch \cite{\refHL}), Theorem 5.1 yields 
$$ 
\|\psi(t)\|_{L^1(\RR)} \leq \|\psi(s)\|_{L^1(\RR)}, 
\qquad 0 \leq s \leq t,
\tag 2.3  
$$
where we neglected the favorable contribution of the Lax discontinuities
appearing in the left-hand side of (2.2). In particular, (2.3) implies that 
the Cauchy problem for (2.1) admits a unique solution 
(in the class of piecewise constant functions at this stage), 
provided $a$ has no 
rarefaction-shock discontinuities.

%_____________________________________________________________________________ 

On the other hand, 
it is clear that the sign of the function $\psi$ is important for the 
sake of deriving the $L^1$ stability of the solutions $\psi$ of (2.1). 
For instance, if $\psi$ has a constant sign for all $(x,t)$, 
then (2.3) holds as an {\it equality\/}
$$
\|\psi(t)\|_{L^1(\RR)} = \|\psi(s)\|_{L^1(\RR)}, 
\qquad 0 \leq s \leq t,  
$$ 
which implies that the Cauchy problem for (2.1) admits at most one solution $\psi$ 
of a given sign. 

%  Let us denote by $\ZZ(\psi)$ the set of points $(x,t) \in \JJ(\psi)$ 
%  where $\psi$ changes sign. 

%______________________________________________________________________________ 

\demo{Proof} 
Denote by $\PP(E)$ the projection of a subset $E$ of the $(x,t)$-plane on the $t$-axis. 
By definition, any piecewise Lipschitz continuous solution $\psi$ is also 
Lipschitz continuous in time with values in $L^1(\RR)$. So,   
it is enough to derive (2.2) for all $t \notin E:=\PP\bigl(\II(a) \cup \II(\psi)\bigr)$. 
The latter is just a finite set. 
The following is valid in each open interval $I$ such that $I \cap E = \emptyset$.  
 
We denote by $x_j(t)$ for $t \in I$ and $j=1,\cdots, m$ the discontinuity lines 
where the function $\psi(.,t)$ changes sign, with the convention that 
$$ 
(-1)^j \, \psi(x,t) \geq 0 \qquad \text{ for } x \in [x_j(t), x_{j+1}(t)]. 
\tag 2.4 
$$
Set $\psi_j^\pm(t) = \psi_\pm(x_j(t),t)$, $\lam_j(t) = \lam^a(x_j(t),t)$, etc. 
Then by using that $\psi$ solves (2.1) we find (for all $t$ in the interval $I$)  
$$ 
\aligned
& \frac{d}{dt}\int_\RR |\psi(x,t)| \, dx
\\
& =
\frac{d}{dt}\sum_{j=1}^m (-1)^ j \, 
\int_{x_j(t)}^{x_{j+1}(t)} \psi(x,t) \, dx \\
& =
\sum_{j=1}^m (-1)^j \, \left(\int_{x_j(t)}^{x_{j+1}(t)} \del_t \psi(x,t) \, dx
+ \lam_{j+1}(t) \, \psi_{j+1}^-(t) - \lam_j(t) \, \psi_j^+(t) \right) \\
& =
\sum_{j=1}^m (-1)^j \,\left(
\int_{x_j(t)}^{x_{j+1}(t)}- \del_x(a(x,t) \, \psi(x,t)) \, dx
+ 
\lam_{j+1}(t) \, \psi_{j+1}^-(t) - \lam_j(t) \, \psi_j^+(t)\right) \\
& =
\sum_{j=1}^m (-1)^j \, \Bigl((a_j^+(t) - \lam_j(t)) \, \psi_j^+(t) +
(a_j^-(t) - \lam_j(t)) \, \psi_j^-(t)\Bigr). 
\endaligned
$$ 
The Rankine-Hugoniot relation associated with (2.1) reads 
$$
(a_j^+(t) - \lam_j(t)) \, \psi_j^+(t)=(a_j^-(t) - \lam_j(t)) \, \psi_j^-(t),
\tag 2.5
$$
therefore by (2.4) 
$$
\frac{d}{dt}\int_\RR |\psi(x,t)| \, dx 
=
2 \, \sum_{j=1}^m \pm (a_j^\pm(t) - \lam_j(t)) \, |\psi_j^\pm(t)|.
\tag 2.6  
$$

Consider each point $x_j(t)$ successively. 
If $x_j(t)$ is a Lax discontinuity, then $a_j^-(t) > \lam_j(t) > a_+(t)$ 
and both coefficients $\pm (a_j^\pm(t) - \lam_j(t))$ are negative. 
If $x_j(t)$ is a rarefaction-shock discontinuity, then 
$a_j^-(t) < \lam_j(t) < a_+(t)$ and the coefficients $\pm (a_j^\pm(t) - \lam_j(t))$ 
are positive. These two cases lead us to the two sums in (2.2). Indeed one just needs 
to observe the following: if $(x, \tau)$ correspond to a Lax or rarefaction-shock 
discontinuity of the speed $a$, but $\psi$ does not change sign at $(x,\tau)$ 
(so it is not counted in (2.6)), then actually by the Rankine-Hugoniot relation 
(see (2.5)) we conclude easily that  
$$
\psi_-(x, \tau) = \psi_+(x, \tau)) = 0, 
$$
and so it does not matter to include the point $(x,\tau)$ in the sums (2.2). 

Suppose next that $x_j(t)$ is an undercompressive discontinuity.
Then the two sides of (2.5) have different sign, therefore 
$$ 
(a_j^+(t) - \lam_j(t)) \, \psi_j^+(t) 
= 
(a_j^-(t) - \lam_j(t))\, \psi_j^-(t) =0, 
$$ 
and the corresponding term in (2.6) vanishes.  $\quad\qed$ 
\enddemo 

%_____________________________________________________________________________ 

Our objective now is to derive an improved version of Theorem 2.2, based 
on a weighted $L^1$ norm adapted to the equation (2.1). 
For piecewise constant functions, we set 
$$
\|\psi(t) \|_{w(t)} := \int_\RR |\psi(x,t)| \, w(x,t) \, dx, 
\tag 2.7
$$
where $w=w(x,t) > 0$ is a piecewise constant and uniformly bounded function.  
We determine this function based on the following constrain on its jumps,  
at each discontinuity of the speed $a$,   
$$ 
w_+(x,t) - w_-(x,t) 
= \cases \leq 0  & \text{ if } (x,t) \in \SS(a), \\
                \\ 
          \geq 0   & \text{ if } (x,t) \in \FF(a). \\
\endcases
\tag 2.8 
$$ 
%%%%%%%% - |a_+(x,t) - a_-(x,t)|
%%%%%%%%   |a_-(x,t) - a_+(x,t)|
The weight is chosen so that the left-hand trace of a slow undercompressive discontinuity 
and the right-hand trace of a fast one are weighted more. This is consistent with 
the immediate observation that 
the terms $\bigl(\lam_j(t) - a_j^-(t) \bigr) \, |\psi_j^-(t)|$ 
and $\bigl(a_j^+(t) - \lam_j(t)\bigr) \, |\psi_j^+(t)|$ 
have a favorable (negative) sign for slow and fast 
undercompressive discontinuities, respectively. On the other hand, 
the jumps of $w$ at Lax or rarefaction-shock discontinuities will remain unconstrained. 
This choice is motivated by the two observations: (i) \, Lax shocks already provide us 
with a good contribution in (2.2), and (ii) \, rarefaction shocks are 
the source of instability and non-uniqueness and cannot be ``fixed up". 

%______________________________________________________________________________ 

The constrain in (2.8) is different for slow and for fast undercompressive 
discontinuities. To actually exhibit a (uniformly bounded) weight
satisfying (2.8), we put a restriction on how the nature of 
the discontinuities change in time as wave interactions take place.  
(An incoming wave may be a slow undercompressive one and become a fast one 
after the interaction, etc. A different constrain is placed before and after 
the interaction.)   

Precisely, we suppose that, to the speed $a=a(x,t)$,
 we can associate on one hand a function 
$\kap : \RR \times \RR_+ \to \RR$ having bounded total variation and 
such that $\JJ(\kap) \subset \JJ(a)$ and $\II(\kap) \subset \II(a)$,  
and on the other hand a partition of the discontinuities  
$$ 
\JJ(a) = \JJ^I(a) \cup \JJ^{II}(a), 
\tag 2.9 
$$ 
so that, for each $(x,t) \in \JJ(a)$, the limits 
$\kap_\pm = \kap_\pm(x,t)$ determine if the wave is slow or fast on its 
left or right side, as follows: 
$$ 
\sgn\bigl(a_\pm(x,t) - \lam(x,t)\bigr) 
= \cases 
  \, \, \, \,  \sgn \kap_\mp & \text{ if } (x,t) \in \JJ^I(a), \\
 - \sgn \kap_\mp  & \text{ if } (x,t) \in \JJ^{II}(a). \\
\endcases
\tag 2.10 
$$ 
Here we use $\sgn (y) = -1, 0, 1$ iff $y < 0, y = 0, y > 0$, respectively. 
Therefore a discontinuity $(x,t) \in \JJ^I(a)$ (for instance) is 
$$
\aligned 
\text{ a Lax one iff }                   & \, \kap_-  <   0 \, \text{ and } \, \kap_+   >  0, \\
\text{ a slow undercompressive one iff } & \, \kap_- \geq 0 \, \text{ and } \, \kap_+ \geq 0, \\
\text{ a fast undercompressive one iff } & \, \kap_- \leq 0 \, \text{ and } \, \kap_+ \leq 0, \\
\text{ a rarefaction-shock iff }         & \, \kap_-   >  0 \, \text{ and } \, \kap_+ < 0.    \\
\endaligned  
$$ 

%_____________________________________________________________________________ 

Furthermore, 
to measure the strength of the jumps, we introduce a piecewise constant 
function, $b=b(x,t)$, having the same jump points as the function $a$. 
For instance, we could assume that there exist constants $C_1, C_2>0$ 
such that at each discontinuity of $a$  
$$
C_1 \, |a_+(y,t) - a_-(y,t)| \leq |b_+(y,t) - b_-(y,t)| \leq C_2 \, |a_+(y,t) - a_-(y,t)|. 
\tag 2.11 
$$
However, strictly speaking, 
this condition will not be used, in the present section at least. 

Based on the functions $\kap$ and $b$ and for $t$ except wave interaction times, 
we can set 
$$
\aligned
V^I(x,t) &   = \sum_{(y,t) \in \JJ^I(a),  \atop y < x} |b_+(y,t) - b_-(y,t)|, 
\\
V^{II}(x,t) & = \sum_{(y,t) \in \JJ^{II}(a), \atop y < x} |b_+(y,t) - b_-(y,t)|,
\endaligned
\tag 2.12  
$$
so that the total variation of $b(t)$ on the interval $(-\infty, x)$ 
decomposes into 
$$
TV_{-\infty}^x(b(t)) = V^I(x,t) + V^{II}(x,t).
\tag 2.13 
$$   
Fix some parameter $m \geq 0$. 
Consider now the weight-function defined for each $(x,t) \in \CC(a)$ by 
$$
w(x,t) = \cases
m + V^I(\infty, t) - V^I(x,t)  
         + V^{II}(x,t) & \text {if } \kap(x,t) >0, \\
\\ 
m + V^I(x,t)
         + V^{II}(\infty, t) - V^{II}(x,t)   & \text {if } \kap(x,t)  \leq 0.\\
\endcases 
\tag 2.14
$$ 
It is immediate to see that indeed (2.8) holds and that with (2.11) 
$$ 
m \leq w(x,t) \leq m + TV(b(t)) \leq m + C_2 \, TV(a(t)), \qquad x \in \RR. 
\tag 2.15 
$$ 
Note also that the weight depends on $b$ and $a$, but not 
on the solution.

\proclaim{Theorem 2.3} Consider a piecewise constant speed $a=a(x,t)$ admitting 
a decomposition $(2.9)$-$(2.10)$ and satisfying the total variation estimate $(2.15)$. 
Consider the weight function $w=w(x,t)$ defined by $(2.13)$. 
Let $\psi$ be any piecewise constant solution of the linear hyperbolic equation $(2.1)$. 
Then the weighted norm $(2.7)$ satisfies for all $0 \leq s \leq t$  
$$ 
\aligned 
& \|\psi(t)\|_{w(t)} 
\\
& +  
\int_s^t \hskip-.3cm  \sum_{(x,\tau) \in \LL(a)} \hskip-.3cm 
  \Big(2 \, m + TV(b) - |b_+(x,\tau) - b_-(x,\tau)| \Big) \, \bigl|a_-(x,\tau) - \lam(x,\tau) \bigr| \,
|\psi_-(x,\tau)| \, d\tau \\ 
& + 
\int_s^t \sum_{(x,\tau) \in \SS(a) \cup \FF(a)} 
|b_+(x,\tau) - b_-(x,\tau)| \, \bigl|a_-(x,\tau) - \lam(x,\tau) \bigr| \,
|\psi_-(x,\tau)| \, d\tau \\ 
& = 
\|\psi(s)\|_{w(s)} + 
\int_s^t \sum_{(x,\tau) \in \RS(a)} 
 \Big(2\, m + TV(b) \Big) \, \bigl|a_-(x,\tau) - \lam(x,\tau) \bigr| \,
|\psi_-(x,\tau)| \, d\tau \\
& \quad +
\int_s^t \sum_{(x,\tau) \in \RS(a)} 
|b_+(x,\tau) - b_-(x,\tau)| \, \bigl|a_-(x,\tau) - \lam(x,\tau) \bigr| \, 
|\psi_-(x,\tau)| \, d\tau. \\
\endaligned 
\tag 2.16
$$
\endproclaim 

The statement (2.16) is sharper than (2.2), as {\it all\/}  
discontinuities contribute now to the decrease of the weighted $L^1$ norm.  
Note that as $m \to \infty$, we recover exactly (2.2) from (2.16).

%_____________________________________________________________________________________________ 

\demo{Proof} We proceed similarly as in the proof of Theorem 2.2.
However, $x_j(t)$ for $t \in I$ (some open interval avoiding the interaction points in $a$ or $\psi$) 
denote now all the jump points in either $a$ or $\psi$. 
We obtain as before the identity 
$$
\aligned 
& \frac{d}{dt} \int_\RR |\psi(x,t)| \, w(x,t) \, dx \\
& = \sum_{j=1}^m 
\Big( (\lam_j(t) - a_j^-(t)) \, |\psi_j^-(t)| \, w_j^-(t)
+ (a_j^+(t) - \lam_j(t)) \, |\psi_j^+(t)|  \, w_j^+(t)  \Big)\\
& = \sum_{j=1}^m 
\Big( \sgn\bigl(\lam_j(t) - a_j^-(t)\bigr) \, w_j^-(t)
    + \sgn\bigl(a_j^+(t) - \lam_j(t)\bigr) \, w_j^+(t) \Big) 
   \, \bigl|\lam_j(t) - a_j^-(t)\bigr| \, |\psi_j^-(t)| \,, 
\endaligned 
\tag 2.17
$$
where we used the Rankine-Hugoniot relation (2.5).

If $x_j(t)$ is a Lax discontinuity in $\JJ^I(a)$, then by (2.11) 
we have $\kap_-<0$ and $\kap_+>0$. So by (2.14) we find 
$$ 
\aligned 
w_j^- & = m + V^I(x_j(t)-) + V^{II}(\infty) - V^{II}(x_j(t)-), \\
w_j^+ & = m + V^I(\infty) - V^I(x_j(t)+) + V^{II}(x_j(t)+),
\endaligned 
$$ 
and so 
$$
\aligned 
&\sgn\bigl(\lam_j(t) - a_j^-(t)\bigr) \, w_j^-(t)
    + \sgn\bigl(a_j^+(t) - \lam_j(t)\bigr) \, w_j^+(t) 
\\
& = 
- w_j^-(t) - w_j^+(t) \\ 
& = 
- 2 \, m - TV(b) + |b_j^+(t) - b_j^-(t)|. 
\endaligned 
\tag 2.18a  
$$

If $x_j(t)$ is a rarefaction-shock discontinuity in $\JJ^I(a)$, then by (2.11) 
we have $\kap_->0$ and $\kap_+<0$. By (2.13) we find 
$$ 
\aligned 
w_j^- & = m + V^I(\infty) - V^I(x_j(t)-) + V^{II}(x_j(t)-), \\
w_j^+ & = m + V^I(x_j(t)+) + V^{II}(\infty) - V^{II}(x_j(t)+),
\endaligned 
$$ 
and so 
$$
\aligned 
&\sgn\bigl(\lam_j(t) - a_j^-(t)\bigr) \, w_j^-(t)
    + \sgn\bigl(a_j^+(t) - \lam_j(t)\bigr) \, w_j^+(t) 
\\
& = w_j^-(t) + w_j^+(t) \\ 
& = 2 \, m + TV(b) + |b_j^+(t) - b_j^-(t)|. 
\endaligned 
\tag 2.18b  
$$

If $x_j(t)$ is a fast undercompressive discontinuity in $\JJ^I(a)$, then by (2.11) 
we have $\kap_- \leq 0$ and $\kap_+ \leq 0$. By (2.13) we find 
$$ 
\aligned 
& \sgn\bigl(\lam_j(t) - a_j^-(t)\bigr) \, w_j^-(t)
    + \sgn\bigl(a_j^+(t) - \lam_j(t)\bigr) \, w_j^+(t) \\
& = w_j^-(t) - w_j^+(t)  \\ 
& = 
m + V^I(x_j(t)-) + V^{II}(\infty) - V^{II}(x_j(t)-) 
\\
& \quad - m - V^I(x_j(t)+) - V^{II}(\infty) + V^{II}(x_j(t)+) \\ 
& =
- |b_j^+(t) - b_j^-(t)|. 
\endaligned 
\tag 2.18c  
$$
Similarly for slow undercompressive discontinuities in $\JJ^I(a)$ we obtain 
$$ 
\sgn\bigl(\lam_j(t) - a_j^-(t)\bigr) \, w_j^-(t)
    + \sgn\bigl(a_j^+(t) - \lam_j(t)\bigr) \, w_j^+(t)
= 
- |b_j^+(t) - b_j^-(t)|. 
\tag 2.18d
$$   

Using (2.18) in (2.17) we conclude that 
$$ 
\aligned 
&\|\psi(t)\|_{w(t)} 
\\
& + 
\int_s^t \sum_{(x,\tau) \in \LL(a)} \Big(2 \, m + TV(b) - |b_+(y,\tau) - b_-(y,\tau)| \Big)  
\, \bigl|a_-(x,\tau) - \lam(x,\tau) \bigr| \, |\psi_-(x,\tau)| \, d\tau \\ 
& + 
\int_s^t \sum_{(x,\tau) \in \SS(a) \cup \FF(a)} 
|b_+(y,\tau) - b_-(y,\tau)| \, \bigl|a_-(x,\tau) - \lam(x,\tau) \bigr| \, |\psi_-(x,\tau)| \, d\tau \\ 
& = 
\|\psi(s)\|_{w(s)} \\
& + 
\int_s^t \sum_{(x,\tau) \in \RS(a)} \Big(2 \, m + TV(b) + |b_+(y,\tau) - b_-(y,\tau)| \Big)  
\, \bigl|a_-(x,\tau) - \lam(x,\tau) \bigr| \, |\psi_-(x,\tau)| \, d\tau, 
\endaligned 
$$
which is equivalent to (2.16). 
$\quad\qed$ 
\enddemo 

Using that $\RS(a)$ is included in the set of points where $\psi$ changes sign, 
it is easy to deduce from (2.16) that:

\proclaim{Corollary 2.4} Under the assumptions and notations in Theorem 2.3, we have 
for all $0 \leq s \leq t$  
$$  
\aligned
& \|\psi(t)\|_{w(t)} 
\\
& \leq 
\|\psi(s)\|_{w(s)} + \Big(2\, m + TV(b) \Big) \, 
\sup_{(x,\tau) \in \RS(a) \atop s \leq \tau \leq t} 
\bigl|b_+(x,\tau) - b_-(x,\tau) \bigr| \, \int_s^t TV(\psi(\tau)) \, d\tau
\endaligned
$$
and, in particular, letting $m \to \infty$   
$$ 
\|\psi(t)\|_{L^1(\RR)} \leq \|\psi(s)\|_{L^1(\RR)} 
+ 2  
\sup_{(x,\tau) \in \RS(a) \atop s \leq \tau \leq t} 
\bigl|b_+(x,\tau) - b_-(x,\tau) \bigr| \, \int_s^t TV(\psi(\tau)) \, d\tau.  
\tag 2.19 
$$
\endproclaim 

Finally, in view of Corollary 2.4, in case the function $a$ contains no rarefaction shocks, 
we deduce that   
$$  
\|\psi(t)\|_{w(t)} \leq \|\psi(s)\|_{w(s)}, \qquad 0 \leq s \leq t.
$$

Observe that this result is achieved, based on a weight that depends on 
an arbitrary function, $b$, and on the sole assumption that a decomposition (2.9)-(2.10) 
of the jumps of $a$ is available. However, our result in this section covers
 only 
piecewise constant solutions. We will see in Section 5 that a stronger structure assumption  
on the coefficients $a$ is necessary to handle general solutions of bounded variation.

%====================================================================================== 
\heading{3. Sharp $L^1$ Estimate for Hyperbolic Conservation Laws}
\endheading

In this section, we apply Theorem 2.3 to the case 
that $a$ is the averaging coefficient (1.5) based on two entropy solutions of (1.1). 
First, we check that the assumptions required in Section 2 on the coefficient 
$a$ do hold in this situation. Therefore Theorem 2.3 
applies to the piecewise constant solutions defined by the wave-front 
traking (also called polygonal approximation) algorithm 
proposed by Dafermos in \cite{\refDafermosone}. 
Next, we observe that, with a suitable choice of the definition of the wave strengths, 
the weighted norm in Section 2 reduces to Liu-Yang's functional. 
Finally we rigorously justify the passage to the limit in the estimate of Theorem 2.3
when the number of wave fronts tends to infinity and exact entropy
solutions of (1.1) are recovered.

Consider the nonlinear scalar conservation law:
$$
\del_t u + \del_x f(u) =0, \qquad u(x,t) \in \RR, 
\tag 3.1
$$
where the flux $f: \RR \to \RR$ is a smooth function. 
Let $u^I$ and $u^{II}$ be two bounded entropy solutions of (3.1) having bounded 
total variation. Given $h>0$ let us approximate the data $u^I(0)$ and $u^{II}(0)$ 
by piecewise constant functions $u^{I,h}(0)$, $u^{II,h}(0)$,  
having finitely many jumps and such that as $h \to 0$  
$$
u^{I,h}(0) \to u^I(0), \quad  u^{II,h}(0) \to u^{II}(0)  
\qquad \text{ in the $L^1$ norm,} 
\tag 3.2  
$$
$$ 
TV(u^{I,h}(0)) \to TV(u^I(0)), \quad  
TV(u^{II,h}(0)) \to TV(u^{II}(0)).  
\tag 3.3   
$$ 

Applying Dafermos' scheme \cite{\refDafermosone}, we can 
construct corresponding, piecewise constant, approximate solutions $u^{I,h}$ and $u^{II,h}$ 
having finitely many jump lines and for $t \geq s \geq 0$ and $p \in [1, \infty]$ 
$$
\|u^{I,h}(t)\|_{L^p(\RR)} \leq \|u^{I,h}(s)\|_{L^p(\RR)},  
\qquad 
\|u^{II,h}(t)\|_{L^p(\RR)} \leq \|u^{II,h}(s)\|_{L^p(\RR)}, 
\tag 3.4  
$$
and for all $-\infty \leq A + M \, (t-s) \leq B - M \, (t-s)$ 
$$
\aligned 
& TV_{A+ M \,(t-s)}^{B-M \, (t-s)} \bigl(u^{I,h}(t)\bigr) \leq 
TV_A^B \bigl(u^{I,h}(s)\bigr), \\ 
& TV_{A+ M \,(t-s)}^{B-M \, (t-s)} \bigl(u^{II,h}(t)\bigr) 
\leq 
TV_A^B \bigl(u^{II,h}(s)\bigr). 
\endaligned 
\tag 3.5
$$ 
More precisely, 
the functions $u^{I,h}$ and $u^{II,h}$ are exact solutions of (3.1) 
satisfying therefore the Rankine-Hugoniot relation at every jump. 
They contain two kinds of jump discontinuities: 
{\it Lax shocks\/} satisfy the so-called Oleinik entropy inequalities, 
while {\it rarefaction jumps\/} do not, but have small strength, that is  
$$
|u^{I,h}(x+,t) - u^{I,h}(x-,t)| \leq h, \qquad 
|u^{II,h}(x+,t) - u^{II,h}(x-,t)| \leq h.   
\tag 3.6 
$$ 
Furthermore, for a subsequence $h \to 0$ at least, we have for each time $t \geq 0$ 
$$
u^{I,h}(t) \to u^I(t), \quad 
u^{II,h}(t) \to u^{II}(t) \qquad \text{ in the $L^1$ norm.}
$$ 
To study the $L^1$ distance between these approximate solutions, we set 
$$
\psi := u^{II,h} - u^{I,h}, 
$$
which is one solution of the linear hyperbolic equation
$$
\del_t \psi + \del_x \bigl( a^h \, \psi\bigr) = 0,  
\qquad 
a^h(x,t) := {f(u^{II,h}(x,t)) - f(u^{I,h}(x,t)) \over u^{II,h}(x,t) - u^{I,h}(x,t)}.  
\tag 3.7 
$$ 

First of all, based on Theorem 2.2 and (3.5)-(3.6), we obtain immediately: 

\proclaim{Theorem 3.1} The approximate solutions $u^{I,h}$ and $u^{II,h}$ satisfy 
the following $L^1$ stability estimate for all $0 \leq s \leq t$  
$$
\aligned 
& \|u^{II,h}(t) - u^{I,h}(t)\|_{L^1(\RR)} 
\\
&
+ 
\int_s^t \sum_{(x,\tau) \in \LL(a)} 
2 \, \bigl(a^h(x-,\tau) - \lam^{a^h}(x,\tau) \bigr) \, |u^{II,h}(x-, \tau) - u^{I,h}(x-, \tau)| \, d\tau \\ 
& \leq  
\|u^{II,h}(s) - u^{I,h}(s)\|_{L^1(\RR)} + 2 h \,  (t-s) \, 
\|f''\|_\infty \, \Big(TV(u^{I,h}(0)) + TV(u^{II,h}(0))\Big).  
\endaligned 
\tag 3.8  
$$
\endproclaim 

{}From the functions $u^I$ and $u^{II}$ we define the function $a$ as in (3.7). 
Recall that the wave front tracking scheme 
converge locally uniformly (see the proof of Theorem 3.5 below for a the 
definition), so that the BV solutions $u^I$ and $u^{II}$ 
are endowed with additional regularity properties. 
Consider for instance the function $u^I$. 
In particular, for all but countably many times $t$ and for each $x$, 
either $x$ is a point of continuity of $u^I$ in the classical sense 
(say $(x,t) \in \CC(u^I)$) or else it is a point of jump in the classical sense 
(say $(x,t) \in \JJ(u^I)$) 
and, to the discontinuity, one can also associate a shock speed, denoted by 
$\lambda^I(x,t)$. 

{}From the properties shared by $u^I$ and $u^{II}$, one deduces immediately 
a similar property for the coefficient $a$. 
Excluding countably many times at most, at each point of jump of $a$ 
we can define the propagation speed $\lam^{a}(x,t)$ of
the discontinuity located at the point $(x,t)$. Namely, we have 
$$
\lam^{a}(x,t) = \cases 
\lambda^I(x,t)      & \text{ if } (x,t) \in \JJ(u^I), \\
\lambda^{II}(x,t)   & \text{ if } (x,t) \in \JJ(u^{II}). \\
\endcases
$$ 
In the limit $h \to 0$ we deduce from $(3.8)$ that: 

\proclaim{Corollary 3.2}  For all $0 \leq s \leq t$ we have 
$$
\aligned 
& \|u^{II}(t) - u^I(t)\|_{L^1(\RR)} 
\\
& + 
\int_s^t \sum_{(x,\tau) \in \LL(a)} 
2 \, \bigl(a(x-,\tau) - \lam^{a}(x,\tau) \bigr) \, |u^{II}(x-, \tau) - u^I(x-, \tau)| \, d\tau \\ 
& \leq  
\|u^{II}(s) - u^I(s)\|_{L^1(\RR)}.   
\endaligned 
\tag 3.9 
$$
\endproclaim 
We omit the proof of Corollary 3.2 as (3.9) is a consequence of  
a stronger estimate proven in Theorem 3.5 below (by taking $m \to \infty$ in (3.15)). 
Note that (3.9) is a stronger statement than the standard $L^1$ contraction estimate
$$ 
\|u^{II}(t) - u^I(t)\|_{L^1(\RR)} 
\leq  
\|u^{II}(s) - u^I(s)\|_{L^1(\RR)}. 
$$ 

\demo{Proof} We apply the estimate (2.2) with $\psi$ replaced with $u^{II,h} - u^{I,h}$. 
We just need to observe (see \cite{\refHL}) 
that all the rarefaction-shock discontinuities in $a^h$ 
are due to rarefaction fronts in $u^{I,h}$ or in $u^{II,h}$, which have small strength 
according to (3.6). In other words we have 
$$ 
\aligned 
& \int_s^t \sum_{(x,\tau) \in \RS(a)} 
2 \, \bigl(\lam^a(x,\tau) - a_-(x,\tau) \bigr) \, |\psi_-(x, \tau)| \, d\tau 
\\
& \leq 
  \sup_{(x,\tau) \in \RS(a) \atop s \leq \tau \leq t} 
 2\,\bigl|a_+(x,\tau) - a_-(x,\tau) \bigr| \, \int_s^t TV(\psi(\tau)) \, d\tau \\ 
& \leq 2\,\|f''\|_\infty \, h \,  \int_s^t TV(\psi(\tau)) \, d\tau \\ 
& \leq 2\,\|f''\|_\infty \, h \,  (t-s) \, \Big(TV(u^{I,h}(0)) + TV(u^{II,h}(0))\Big). 
\endaligned 
$$ 
This establishes (3.8). 
$\quad\qed$ 
\enddemo

%_____________________________________________________________________________ 

We now want to apply Theorem 2.3 and control a weighted norm of $u^{II,h} - u^{I,h}$. 
In this direction, our main observation is: 

\proclaim{Lemma 3.3} When the function $f$ is strictly convex, 
the coefficient $a^h$ satisfies all of the assumptions $(2.9)$-$(2.10)$.
\endproclaim

\demo{Proof} The function $a^h$ is piecewise constant, and we can associate to this function 
an obvious decomposition of the form (2.9). To establish 
(2.10), consider for instance a jump point $(x,t) \in \JJ(u^{I,h}) \cap \CC(u^{II,h})$, 
together with its left- and right-hand traces $u^I_-$ and $u^I_+$. Since $u^{I,h}$ 
is a solution of (3.1), 
the corresponding speed $\lam = \lam (x,t)$ satisfies the Rankine-Hugoniot relation: 
$$
- \lam \, \bigl(u^I_+ - u^I_- \bigr) + f(u^I_+)  - f(u^I_-)  = 0.  
$$
Thus the term in the left-hand side of (2.10) takes the form
$$ 
\aligned 
a_\pm(x,t) - \lam(x,t) 
& = {f(u^{II}) - f(u_\pm^I) \over u^{II} - u_\pm^I} 
- {f(u_+^I) - f(u_-^I) \over u_+^I - u_-^I} \\
& = \int_0^1 \Big( 
    f'\bigl(\theta \, u^{II} + (1 - \theta) \, u_\pm^I\bigr) 
    - 
    f'\bigl(\theta \, u_\mp^I + (1 - \theta) \, u_\pm^I\bigr) 
\Big)  \, d\theta. 
\endaligned 
$$ 
Thus we obtain 
$$
\aligned 
& a_\pm(x,t) - \lam(x,t) 
= \mu \, \bigl(u^{II} - u_\mp^I \bigr), \\  
& \mu : = \int_0^1 \int_0^1 
    f''\Big(\rho \, \bigl(\theta \, u^{II} + (1 - \theta) \, u_\pm^I\bigr) 
+ (1 - \rho) \, \bigl(\theta \, u_\mp^I + (1 - \theta) \, u_\pm^I\bigr) \Big)
\, \theta \, d\theta d\rho. 
\endaligned 
\tag 3.10 
$$ 
Since $f$ is strictly convex, the coefficient is bounded away from zero. 
In view of (3.10), if we now choose $\kap(x,t) := u^{II,h} - u^{I,h}$,
the desired property (2.10) holds true. 
$\quad\qed$ 
\enddemo 

%_________________________________________________________________________________ 

Next, we define the weight $w^h= w^h(x,t)$ associated with the function $a^h$, 
by the formula (2.14) in which we specify 
$$
\kap^h(x,t) := u^{II,h} - u^{I,h}. 
\tag 3.11
$$
It follows immediately from Theorem 2.3 that:

\proclaim{Theorem 3.4} Suppose that the function $f$ is strictly convex. 
The approximate solutions constructed by Dafermos scheme 
satisfy the $L^1$ stability estimate for all $0 \leq s \leq t$  
$$ 
\aligned 
& \|u^{II,h}(t) - u^{I,h}(t)\|_{w^h(t)} \\
& + 
\int_s^t \sum_{(x,\tau) \in \LL(a^h)} 
  \Big(2 \, m + TV(b^h) - |b^h(x+,\tau) - b^h(x-,\tau)| \Big) 
  \\
  & \hskip3.cm \, \bigl|a^h(x-,\tau) - \lam^h(x,\tau) \bigr| 
   \, |u^{II,h}(x-, \tau) - u^{I,h}(x-, \tau)| \, d\tau \\ 
& + 
\int_s^t \sum_{(x,\tau) \in \SS(a^h)\cup\FF(a^h)} 
|b^h(x+,\tau) - b^h(x-,\tau)| \, \bigl|a^h(x-,\tau) - \lam^h(x,\tau) \bigr| \, 
   \\
  & \hskip3.cm    |u^{II,h}(x-, \tau) - u^{I,h}(x-, \tau)| \, d\tau \\ 
& = 
\|u^{II,h}(s) - u^{I,h}(s)\|_{w^h(s)} \\
& + 
\int_s^t \sum_{(x,\tau) \in \RS(a^h)} 
 \Big(2\, m + TV(b^h) + 
|b^h(x+,\tau) - b^h(x-,\tau)|\Big)   \\
  & \hskip3.cm \, \bigl|a^h(x-,\tau) - \lam^h(x,\tau) \bigr| \, 
   |u^{II,h}(x-, \tau) - u^{I,h}(x-, \tau)|  \, d\tau, 
\endaligned 
\tag 3.12 
$$ 
where $a^h$ is the averaging coefficient defined in $(3.7)$ and
$\lam^h(x,\tau)$ represents 
the speed of the discontinuity located at $(x,\tau) \in \JJ(a^h)$. 
\endproclaim

We emphasize that (3.12) is an {\it equality\/} in which the contribution to the $L^1$ norm 
of each type of wave appears clearly. The coefficient $a^h$ exhibits three types of waves: 
the Lax and undercompressive discontinuities in $a^h$ contribute to the decay 
of the $L^1$ weighted distance. The statement (3.12) quantifies sharply this effect. 
On the other hand, the rarefaction-shocks appearing in the right-hand side of (3.12) increase
the $L^1$ norm. 

In the rest of this section, we assume that the function $b=b^h$ is chosen 
to be specifically  
$$ 
b^h(x+,t) - b^h(x-,t) = \cases 
u^{I,h}(x+,t)    - u^{I,h}(x-,t)         & \text{ if } (x,t) \in \JJ(u^{I,h}), \\ 
u^{II,h}(x+,t) - u^{II,h}(x-,t)      & \text{ if } (x,t) \in \JJ(u^{II,h}), \\ 
\endcases
\tag 3.13 
$$ 
but a more general definition is possible. 

%_______________________________________________________________________

Our next purpose is to pass to the limit ($h \to 0$) in the statement
established in Theorem 3.4 for piecewise constant 
approximate solutions. We recover here a result derived by 
Dafermos \cite{\refDafermosthree} via a different approach. 
Recall the notation $\CC(u^I)$, $\SS(u^I)$, etc introduced earlier. 
Denote by $\II(u^I)$ the countable set of interactions times. 
Let $V^I(t)$ be the total variation function associated with $u^I(t)$. 
Based on the functions $V^I(t)$ and $V^{II}(t)$,
we then define the weight $w$ as in (2.14) but with (2.12) replaced by the 
total variation functions of $u^I(t)$ and $u^{II}(t)$, with 
$\kap := u^{II} - u^I$ and 
$$ 
b(x+,t) - b(x-,t) = \cases 
u^I(x+,t)    - u^I(x-,t)         & \text{ if } (x,t) \in \JJ(u^I), \\ 
u^{II}(x+,t) - u^{II}(x-,t)      & \text{ if } (x,t) \in \JJ(u^{II}). \\ 
\endcases
\tag 3.14 
$$ 

Furthermore, to any functions of bounded variation $u,v,w$ 
in the space variable $x$
(the time variable being fixed) we associate the measure on $\RR$ 
$$ 
\mu = \bigl(a(u,v) - f'(u) \bigr) \, (v - u) \, dw  
$$ 
understood as the nonconservative product in the sense 
of Dal Maso, LeFloch and Murat \cite{\refDLM} and characterized by the 
following two conditions: 
\roster
\item If $B$ is a Borel set included in the set of continuity points of $w$
$$
\mu (B) = \int_B  \bigl(a(u,v) - f'(u) \bigr) \, (v - u) \, dw, 
\tag 3.15a
$$
where the integral is defined in a classical sense; 
\item If $x$ is a point of jump of $w$, then 
$$
\aligned
\mu(\{x\}) = {1 \over 2} \, \Bigl(
 &  \bigl(a(u_+, v_+) - a(u_-,u_+) \bigr) \, (v_+ - u_+) 
 \\
 & + 
   \bigl(a(u_-, v_-) - a(u_-,u_+) \bigr) \, (v_- - u_-)  \Bigr) 
\, |w_+ - w_-|  
\endaligned
\tag 3.15b
$$ 
with $u_\pm = u(x\pm)$, etc. 
\endroster 
Note that, if $u=u^I$ and $v=u^{II}$, 
the two terms $\bigl(a(u_\pm, v_\pm) - a(u_-,u_+) \bigr) \, (v_\pm - u_\pm)$
in fact coincide. 

%_______________________________________________________________________

\proclaim{Theorem 3.5} Let the function $f$ be strictly convex and let 
$u^I$ and $u^{II}$ be two entropy solutions of bounded variation of 
the conservation law $(1.1)$. For all $0 \leq s \leq t$ we have    
$$ 
\aligned 
& \|u^{II}(t) - u^I(t)\|_{w(t)} \\
& + 
\int_s^t \hskip-.22cm 
\sum_{(x,\tau) \in \LL(a)\cap\JJ(u^I)} \hskip-.22cm  q \, 
\bigl|a\bigl(u^I(x-), u^{II}(x-)\bigr) - 
a\bigl(u^I(x+), u^I(x-)\bigr) \bigr| \, 
|u^{II}(x) - u^I(x)| \, d\tau \\ 
& + 
\int_s^t \hskip-.22cm  \sum_{(x,\tau) \in \LL(a)\cap\JJ(u^{II})} \hskip-.22cm  q \, 
\bigl|a\bigl(u^I(x-), u^{II}(x-)\bigr) - 
a\bigl(u^{II}(x+), u^{II}(x-)\bigr) \bigr| \, 
|u^{II}(x) - u^I(x)| \, d\tau \\ 
& + 
\int_s^t \int_\RR \bigl(a(u^I, u^{II}) - f'(u^I) \bigr) 
\, (u^{II} - u^I) \, dV^I \, d\tau
\\
&
+ \int_s^t \int_\RR 
 \bigl(a(u^I, u^{II}) - f'(u^{II}) \bigr) \, (u^{I} - u^{II}) 
\, dV^{II} \, d\tau \\
& 
\leq \|u^{II}(s) - u^I(s)\|_{w(s)}. 
\endaligned 
\tag 3.16 
$$  
where 
$q=q(\tau) = 2 \, m + TV(u^I(\tau)) + TV(u^{II}(\tau))$. 
\endproclaim 
%__________________________________________________________________________ 

Observe that the terms in integrals in (3.16) globally contribute to the 
decrease of weighted norm, as is better seen rewriting the formula
as follows ($V_c^I$ and $V_c^{II}$ being the continuous parts of the measures $V^I$ and $V^{II}$): 
$$ 
\aligned 
& \|u^{II}(t) - u^I(t)\|_{w(t)} \\
& + 
\int_s^t \sum_{(x,\tau) \in \LL(a)\cap\JJ(u^I)} 
\bigl(q - |u^I_+ - u^I_-|\bigr) \, 
\bigl|a\bigl(u^I_-, u^{II}_+\bigr) - 
a\bigl(u^I_+, u^I_-\bigr) \bigr| \, 
|u^{II} - u^I| \, d\tau \\ 
& + 
\int_s^t \sum_{(x,\tau) \in \LL(a)\cap\JJ(u^{II})} 
\bigl(q - |u^{II}_+ - u^{II}_-|\bigr) \, 
\bigl|a\bigl(u^I_-, u^{II}_-\bigr) - 
a\bigl(u^{II}_+, u^{II}_-\bigr) \bigr| \, 
|u^{II} - u^I| \, d\tau \\  
& + 
\int_s^t \sum_{(x,\tau) \in \bigl(\SS(a)\cup\FF(a)\bigr)\cap\JJ(u^I)} 
\bigl|a\bigl(u^I_-, u^{II}_-\bigr) - 
a\bigl(u^I_+, u^I_-\bigr) \bigr| \, 
|u^{II} - u^I| \, 
|u^I_+ - u^I_-|
d\tau \\ 
& + 
\int_s^t \sum_{(x,\tau) \in \bigl(\SS(a)\cup\FF(a)\bigr)\cap\JJ(u^{II})} 
\bigl|a\bigl(u^I_-, u^{II}_-\bigr) - 
a\bigl(u^{II}_+, u^{II}_-\bigr) \bigr| \, 
|u^{II} - u^I| \, 
|u^{II}_+ - u^{II}_-|
d\tau \\ 
& +
\int_s^t \int_\RR \bigl|a(u^I, u^{II}) - f'(u^I) \bigr| 
\, |u^{II} - u^I| \, dV^I_c \, d\tau 
\\ &
+ \int_s^t \int_\RR 
 \bigl|a(u^I, u^{II}) - f'(u^{II}) \bigr| \, |u^{I} - u^{II}| 
\, dV^{II}_c \, d\tau \\
& 
\leq \|u^{II}(s) - u^I(s)\|_{w(s)}. 
\endaligned 
\tag 3.16' 
$$

%___________________________________________________________________________
The following estimate is a direct consequence of the definition (3.15):

\proclaim{Lemma 3.6} There exists a constant $C>0$ such that 
for all functions of bounded variation $u, \tilde u, v, \tilde v, w$ defined on some interval $[\alpha, \beta]$   
$$
\aligned 
& \left| \int_\alpha^\beta \bigl(a(u,v) - f'(u) \bigr) \, (v - u) \, dw
- 
 \int_\alpha^\beta 
\bigl(a(\tilde u, \tilde v) - f'(\tilde u) \bigr) \, (\tilde v - \tilde u) \, dw
\right| \\
& \leq 
C \, \left( \|\tilde u - u\|_{L^\infty(\alpha, \beta)} + 
       \|\tilde v - v \|_{L^\infty(\alpha, \beta)} \right) \, TV_{[\alpha,\beta]}(w).  
\endaligned 
\tag 3.17
$$
\endproclaim

%____________________________________________________________________________

\demo{Proof of Theorem 3.5} 

{\bf Step 1 : Preliminaries.} 

For each $t \geq 0$, the functions $V^{I,h}(t)$ and $V^{II,h}(t)$ 
associated with the wave front tracking approximations $u^{I,h}(t)$ 
and $u^{II,h}(t)$ are of uniformly bounded variation as $h \to 0$. 
The measures $dV^{I,h}$ and $dV^{II,h}$ are also Lipschitz continuous in time 
(with constant independent of $h$) for the weak convergence, except at 
interaction points. On the other hand, interaction times in 
the limiting solutions are at most countable. 
Therefore, extracting subsequences if necessary,
the measures $dV^{I,h}$ and $dV^{II,h}$
converge to some limiting (non-negative) measures, say:
$$
dV^{I,h}(t) \to d\bar V^I(t) \qquad 
dV^{II,h}(t) \to d\bar V^{II}(t). 
\tag 3.18  
$$   
By lower semi-continuity, we have at each time $t$ 
$$ 
dV^I(t)    \leq d\bar V^I(t),   
\qquad 
dV^{II}(t) \leq d\bar V^{II}(t), 
\tag 3.19 
$$ 
and, in particular, at each $(x,t)$  
$$ 
V^I(x,t)    \leq \bar V^I(x,t), \qquad 
V^{II}(x,t) \leq \bar V^{II}(x,t).  
\tag 3.20a 
$$
$$ 
\aligned 
& V^I(+\infty,t) - V^I(x,t)   \leq \bar V^I(+\infty,t) - \bar V^I(x,t), \\
& V^{II}(+\infty,t)- V^{II}(x,t) \leq \bar V^{II}(+\infty,t) - \bar V^{II}(x,t).  
\endaligned 
\tag 3.20b 
$$

Based on the functions $\bar V^I(t)$ and $\bar V^{II}(t)$, 
on the coefficient $\kappa:= u^{II} - u^I$ and on the function in (3.14), 
we can define a weight denoted by $\bar w$, along the same lines as in (2.14). 
We will show that the left-hand side of (3.16) is bounded above by 
$$ 
\aligned 
& \|u^{II}(t) - u^I(t)\|_{\bar w(t)} \\
& + 
\int_s^t \hskip-.22cm \sum_{(x,\tau) \in \LL(a)\cap\JJ(u^I)} \hskip-.22cm 
\bar q \, 
\bigl|a\bigl(u^I(x-), u^{II}(x-)\bigr) - a\bigl(u^I(x+), u^I(x-)\bigr) \bigr| \, 
|u^{II}(x) - u^I(x)| \, d\tau \\ 
& + 
\int_s^t \hskip-.22cm \sum_{(x,\tau) \in \LL(a)\cap\JJ(u^{II})} \hskip-.22cm 
\bar q \, 
\bigl|a\bigl(u^I(x-),u^{II}(x-)\bigr)-a\bigl(u^{II}(x+),u^{II}(x-)\bigr)\bigr| \,|u^{II}(x) - u^I(x)| \, d\tau \\ 
& + 
\int_s^t \int_\RR \bigl(a\bigr(u^I, u^{II}\bigr) - f'(u^I)\bigr) \,
 (u^{II} - u^I) \, d\bar V^I(y, \tau) d\tau \\ 
& +
\int_s^t \int_\RR \bigl(a\bigl(u^I, u^{II}\bigr) - f'(u^{II})\bigr) 
\, (u^I - u^{II}) \, d\bar V^{II}(y, \tau) d\tau \\ 
\endaligned 
\tag 3.21
$$
where $\bar q : = 2 \, m + \bar V^I(+\infty) + \bar V^{II}(+\infty)$, 
and that (3.21) coincides 
with the desired upper bound $\|u^{II}(s) - u^I(s)\|_{w(s)}$. 
The former statement is postponed to Step 5 below and we focus now on the latter. 

%_________________________________________________________________________ 

Fix some $t \geq s \geq 0$ and rewrite (3.12) in the equivalent form
$$ 
\aligned 
& \|u^{II,h}(t) - u^{I,h}(t)\|_{w^h(t)} \\
& + 
\int_s^t \sum_{(x,\tau) \in \LL(a^h)} 
  \Big(2 \, m + TV(b^h) \Big) \bigl|a^h(x-,\tau) - \lam^h(x,\tau) \bigr| 
   \, |u^{II,h}(x-, \tau) - u^{I,h}(x-, \tau)| \, d\tau \\ 
& + 
\int_s^t \sum_{(x,\tau) \in \JJ^I(a^h)} 
|b^h(x+,\tau) - b^h(x-,\tau)| \, \bigl(a^h(x-,\tau) - \lam^h(x,\tau) \bigr) \, 
    (u^{II,h}(x-, \tau) - u^{I,h}(x-, \tau)) \, d\tau \\ 
& + 
\int_s^t \sum_{(x,\tau) \in \JJ^{II}(a^h)} 
|b^h(x+,\tau) - b^h(x-,\tau)| \, \bigl(a^h(x-,\tau) - \lam^h(x,\tau) \bigr) \, 
    (u^{I,h}(x-, \tau) - u^{II,h}(x-, \tau)) \, d\tau \\ 
& = 
\|u^{II,h}(s) - u^{I,h}(s)\|_{w^h(s)} \\
& + 
\int_s^t \sum_{(x,\tau) \in \RS(a^h)} 
 \Big(2\, m + TV(b^h) \Big) \, \bigl|a^h(x-,\tau) - \lam^h(x,\tau) \bigr| \, 
   |u^{II,h}(x-, \tau) - u^{I,h}(x-, \tau)|  \, d\tau, \\
\endaligned 
\tag 3.22 
$$ 
or, with obvious notations, 
$$ 
\|u^{II,h}(t) - u^{I,h}(t)\|_{w^h(t)} 
+ \Omega_1^h +  \Omega_2^h
= 
\|u^{II,h}(s) - u^{I,h}(s)\|_{w^h(s)} 
+ \Omega_3^h. 
\tag 3.23   
$$  
As the maximum strength of rarefaction fronts in $u^{I,h}$ and 
$u^{II,h}$ vanishes with $h$ (see (3.6)) 
and rarefaction shocks in $a^h$ arise only from these rarefaction fronts 
(see (1.6)), 
we have 
$$
\Omega_3^h \to 0 \qquad \text{ as } h \to 0. 
\tag 3.24 
$$
On the other hand, we can always choose the (initial) approximations 
at time $s$ in such a way that 
$$
\bar w(s) = w(s)
\tag 3.25
$$
and 
$$ 
\lim_{h \to 0} \|u^{II,h}(s) - u^{I,h}(s)\|_{w^h(s)}  
= \|u^{II}(s) - u^I(s)\|_{w(s)}. 
\tag 3.26  
$$ 
It remains to prove that the limit of the left-hand side 
of (3.22) is exactly (3.21). This will be established in the following 
three steps. 

%________________________________________________________________________

\vskip.5cm

{\bf Step 2 : } We will 
rely on the local uniform convergence of the front tracking approximations 
(see Bressan and LeFloch \cite{\refBL}). 
For all but countably many times $\tau$ we have the following properties for 
$u^I$ (as well as for $u^{II}$): 
\roster
\item For each point of jump $z$ of $u^I$ 
there exists a sequence $z^h \to z$ such that for each $\eps>0$ there exists 
$\delta>0$ such that 
$$
\aligned 
& |u^{I,h} (x) - u^I(z+) | + |u^I (x) - u^I(z+)|  < \eps 
   \quad \text{ for all } x-z^h \in (0, \delta), \\
& |u^{I,h} (x) - u^I(z-) | + |u^I (x) - u^I(z-) | < \eps 
   \quad \text{ for all } x-z^h \in (-\delta, 0)  
\endaligned 
\tag 3.27a
$$ 
and (clearly) 
$$
\aligned 
& |u^I (x) - u^I(z+)|  < \eps 
   \quad \text{ for all } x-z \in (0, \delta), \\
& |u^I (x) - u^I(z-) | < \eps 
   \quad \text{ for all } x-z \in (-\delta, 0). 
\endaligned 
\tag 3.27b 
$$ 
\item For each point of continuity $z$ of $u^I$ and 
for each $\eps>0$, there exists $\delta>0$ such that 
$$ 
|u^{I,h} (x) - u^I(z) | + |u^I (x) - u^I(z)|  < \eps 
\quad \text{ for all } x-z \in (-\delta, \delta). 
\tag 3.28   
$$ 
\endroster 

We also recall from \cite{\refBL} that, for all but countably many times $t$, 
the atomic parts of the measures $\bar V^I$ and $\bar V^{II}$ coincide 
with the one of $V^I$ and $V^{II}$, that is for each $y \in \RR$ 
$$
\aligned 
& \bar V^I(y+, t) - \bar V^I(y-, t) =  V^I(y+, t) - V^I(y-, t), \\
& \bar V^{II}(y+, t) - \bar V^{II}(y-, t) = V^{II}(y+, t) - V^{II}(y-, t). 
\endaligned 
\tag 3.29
$$

Following LeFloch and Liu \cite{\refLL}  
who established the weak stability of nonconservative products under
local uniform convergence, we want to show that  
$$ 
\aligned 
\Omega_2^h(\tau)  
: =& \int_\RR \bigl(a\bigl(u^{I,h}(y, \tau), u^{II,h}(y,\tau)\bigr)
 - f'(u^{I,h}(y,\tau))\bigr) 
\, (u^{II,h}(y, \tau) - u^{I,h}(y, \tau)) \, dV^{I,h}(y)  \\
&+ \int_\RR \bigl(a\bigl(u^{I,h}(y, \tau), u^{II,h}(y,\tau)\bigr) 
- f'(u^{II,h}(y,\tau))\bigr) 
\, (u^{I,h}(y, \tau) - u^{II,h}(y, \tau)) \, dV^{II,h}(y)  \\
\longrightarrow
 &  \int_\RR \bigl(a\bigl(u^I(y, \tau), u^{II}(y,\tau)\bigr) - f'(u^I(y,\tau))\bigr) \, 
(u^{II}(y,\tau) - u^I(y, \tau)) \, d\bar V^I(y) \\
&+\int_\RR \bigl(a\bigl(u^I(y, \tau), u^{II}(y,\tau)\bigr) - f'(u^{II}(y,\tau))\bigr) \, 
(u^{I}(y, \tau) - u^{II}(y, \tau)) \, d\bar V^{II}(y). \\
\endaligned 
\tag 3.30
$$ 
By Lebesgue dominated convergence theorem and 
since a uniform bound in $\tau$ and $h$ is 
available, it will follow from (3.29) that 
$$
\aligned 
\Omega_2^h 
= \int_s^t \Omega_2^h(\tau) \, d\tau  
& \longrightarrow  
  \int_s^t \int_\RR \bigl(a\bigl(u^I, u^{II}\bigr) - f'(u^I)\bigr) \, 
  (u^{II} - u^I) \, d\bar V^I(y, \tau) \, d\tau \\
& + \int_s^t \int_\RR \bigl( a\bigl(u^I, u^{II}\bigr) - f'(u^{II})\bigr) \, 
  (u^{I} - u^{II}) \, d\bar V^{II}(y, \tau) \, d\tau. \\
\endaligned 
\tag 3.30'  
$$

%_____________________________________________________________________ 

Given $\eps>0$, select finitely many (large) jumps in $u^I$ or $u^{II}$, 
located at $y_1, y_2, \ldots y_n$, so that 
$$
\sum_{x \neq y_j \atop 
j=1,2, \ldots, n} \bigl|u^I(x+) - u^I(x-)\bigr|  
+ \bigl|u^{II}(x+) - u^{II}(x-)\bigr| 
 < \eps. 
\tag 3.31 
$$
To each $y_j$ we associate the corresponding discontinuity point $y_j^h$ in $u^{I,h}$
or $u^{II,h}$. 
To simplify the presentation we will focus on the case where 
$y_j<y_j^h<y_{j+1}<y_{j+1}^h$ for all $j$. 
The other cases can be treated similarly.
In view of the local convergence property (3.27)--(3.28) and by extracting 
a covering of the interval $[y_0, y_n]$, we have also  
$$ 
|u^{I,h}(x) - u^I(x)| + |u^{II,h}(x) - u^{II}(x)| \leq 2 \, \eps, 
\quad x \in \bigl(y_j^h, y_{j+1}\bigr) \subseteq \bigl(y_j, y_{j+1}\bigr). 
\tag 3.32   
$$

In view of (3.30) we can construct functions $u^I_\eps$ and $u^{II}_\eps$ 
that are continuous everywhere except possibly at the points $y_j$ 
and such that the following conditions hold 
with $u$ replaced by either $u^I$ or $u^{II}$:
$$
\aligned 
& TV\left(u_\eps; \RR \setminus\bigl\{y_1, \ldots, y_n\bigr\}\right) \leq C \, 
TV\left(u; \RR \setminus\bigl\{y_1, \ldots, y_n\bigr\}\right), \\
& \|u-u_\eps\|_\infty \leq C \, \eps, 
\quad 
TV\left(u-u_\eps; \RR \setminus\bigl\{y_1, \ldots, y_n\bigr\}\right) \leq C \, \eps, 
\endaligned 
\tag 3.33 
$$
where $C$ is independent of $\eps$.

Consider the decompositions 
$$
\int_\RR \bigl(a\bigl(u^{I,h}, u^{II,h}\bigr) - f'(u^{I,h})\bigr) \, 
 (u^{II,h} - u^{I,h}) \, d V^{I,h} 
= 
\sum_{j=0}^n \int_{(y_j^h, y_{j+1}^h)} \cdots 
+ \sum_{j=1}^n \int_{\{y_j^h\}} \cdots
$$
and 
$$
\int_\RR \bigl( a\bigl(u^I, u^{II}\bigr) - f'(u^I)\bigr) 
\, (u^{II} - u^I) \, d\bar V^I
=
\sum_{j=0}^n \int_{(y_j, y_{j+1})} \cdots 
+ \sum_{j=1}^n \int_{\{y_j\}} \cdots. 
$$
Here $y_0^h = y_0 =-\infty$ and $y_{n+1}^h = y_{n+1} = +\infty$. 
Thus in (3.30) we have to estimate 
$$
\aligned 
\Omega_2^h(\tau) 
& = \int_\RR \bigl(a\bigl(u^{I,h}, u^{II,h}\bigr) - f'(u^{I,h})\bigr) \, 
    ( u^{II,h} - u^{I,h}) \, d V^{I,h} 
    \\
    & \quad  -
    \int_\RR \bigl( a\bigl(u^I, u^{II}\bigr)
    - f'(u^I)\bigr) \, (u^{II} - u^I) \, d\bar V^I\\
& = T_1^h + T_2^h 
\endaligned 
\tag 3.34
$$ 
with 
$$
\aligned 
T_1^h := & 
\sum_{j=1}^n \int_{\{y_j^h\}} 
\bigl(a\bigl(u^{I,h}, u^{II,h}\bigr) - f'(u^{I,h})\bigr) \, 
 (u^{II,h} - u^{I,h}) \, d V^{I,h}\\ 
& - 
\sum_{j=1}^n \int_{\{y_j\}}  
\bigl(a\bigl(u^I, u^{II}\bigr) - f'(u^I)\bigr) 
\, (u^{II} - u^I) \, d\bar V^I
\endaligned 
$$
and 
$$ 
\aligned 
T_2^h := &  
\sum_{j=0}^n \int_{(y_j^h, y_{j+1}^h)}
\bigl(a\bigl(u^I, u^{II}\bigr) - f'(u^{I,h})\bigr) \, 
 (u^{II,h} - u^{I,h}) \, d V^{I,h}\\
& -
\sum_{j=0}^n \int_{(y_j, y_{j+1})} 
\bigl(a\bigl(u^I, u^{II}\bigr) - f'(u^I)\bigr) 
\, (u^{II} - u^I) \, d\bar V^I. 
\endaligned 
$$

First, relying on the convergence property (3.29) we have immediately 
$$ 
\aligned 
T_1^h 
& = \sum_{j=1}^n 
 \bigl(a\bigl(u^{I,h}(y^h_j-), u^{II}(y^h_j-)\bigr)
       - \lambda^{I,h}(y^h_j-) \bigr) \, \bigl(u^{II,h}(y^h_j-) - u^{I,h}(y^h_j-)\bigr) 
\, |u^{I,h}(y^h_j+) - u^{I,h}(y^h_j-)| 
\\ 
& - \bigl(a\bigl(u^{I}(y_j-), u^{II}(y_j-)\bigr) 
- \lambda^I(y_j-) \bigr) \, \bigl(u^{II}(y_j-) - u^I(y_j-)\bigr) 
\, |u^I(y_j+) - u^I(y_j-)|,  
\endaligned 
$$
so that 
$$ 
\bigl|T_1^h \bigr| \leq  
C \, \sum_{j=1}^n\sum_\pm |u^{I,h}(y^h_j\pm) - u^I(y_j\pm)| 
 + |u^{II,h}(y^h_j\pm) - u^{II}(y_j\pm)|. 
$$  
Thus, in view of the local convergence at jump points (3.27a), 
for $h$ small enough we obtain   
$$
\bigl|T_1^h\bigr| \leq C \, \eps.
\tag 3.35  
$$

%________________________________________________________ 

Relying on the simplifying 
assumption $y_j<y_j^h<y_{j+1}<y_{j+1}^h$ for all $j$, 
we can decompose $T_2^h$ as follows: 
$$
\aligned 
T_2^h  
& =\sum_{j=0}^n \int_{(y_j^h, y_{j+1})}
  \bigl( a\bigl(u^{I,h}, u^{II,h}\bigr) - f'(u^{I,h})\bigr) \, 
    (u^{II,h} - u^{I,h}) \, d V^{I,h} 
    \\
    & \qquad \qquad - \bigl( a\bigl(u^{I}, u^{II}\bigr) 
- f'(u^I)\bigr) \, (u^{II} - u^I) \, d\bar V^I \\ 
& - 
\sum_{j=0}^n \int_{(y_j, y_j^h]} 
\bigl(a\bigl(u^I, u^{II}\bigr) - f'(u^I)\bigr) \, 
(u^{II} - u^I) \, d\bar V^I \\ 
& + 
\sum_{j=0}^n \int_{[y_{j+1}, y_{j+1}^h)}
\bigl(a\bigl(u^{I,h}, u^{II,h}\bigr) - f'(u^{I,h})\bigr) \, 
 (u^{II,h} - u^{I,h}) \, d V^{I,h}\\ 
& =: T^h_{2,1} + T^h_{2,2} + T^h_{2,3}. 
\endaligned  
\tag 3.36 
$$ 
We first consider $T_{2,2}^h$: 
$$
\aligned 
T_{2,2}^h 
=   
- \sum_{j=0}^n \int_{(y_j, y_j^h]} 
  & \bigl(a\bigl(u^I(y_j+), u^{II}(y_j+)\bigr) - f'(u^I(y_j+))\bigr) 
   \, \bigl(u^{II}(y_j+) - u^I(y_j+)\bigr) \, d\bar V^I(y) \\  
+ \sum_{j=0}^n \int_{(y_j, y_j^h]} 
  & \Big\{
  \bigl(a\bigl(u^I(y), u^{II}(y)\bigr) - f'(u^I(y))\bigr) \, \bigl(u^{II}(y) - u^I(y)\bigr) \\ 
  &  -  \bigl(a\bigl(u^I(y_j+), u^{II}(y_j+)\bigr) - f'(u^I(y_j+))\bigr) \, 
\bigl(u^{II}(y_j+) - u^I(y_j+)\bigr) \Big\} \, d\bar V^I(y).  
\endaligned 
$$
Therefore, with (3.17), we obtain 
$$
\aligned 
\bigl|T_{2,2}^h \bigr| 
\leq & C \, \sum_j \bigl|\bar V^I(y_j+) - \bar V^I(y_j^h+)\bigr| \\
     & + C \,  V^I(+\infty) \, \Big(
       \sup_{y \in (y_j, y_j^h]} |u^I(y) - u^I(y_j+)| 
       + \sup_{x \in (y_j, y_j^h]} |u^{II}(y) - u^{II}(y_j+)| \Big). 
\endaligned 
$$
Since $y_j^h \to y_j$, we have 
$\bigl|\bar V^I(y_j+) - \bar V^I(y_j^h+)\bigr| \to 0$, therefore 
for $h$ sufficiently small 
$$ 
\bigl|T_{2,2}^h\bigr| \leq C \, \eps. 
\tag 3.37  
$$
A similar argument for $T_{2,3}^h$ shows that 
$$
\bigl|T_{2,3}^h\bigr| \leq C \, \eps. 
\tag 3.38   
$$ 
%___________________________________________________________________________ 

Next consider the decomposition 
$$ 
\aligned 
\bigl(a(& u^{I,h}, u^{II,h}) - f'(u^{I,h})\bigr) \, 
  (u^{II,h} - u^{I,h}) \, d V^{I,h}
 - \bigl(a(u^I, u^{II}) - f'(u^I)\bigr) \, (u^{II} - u^I) \, d\bar V^I \\
= \bigl(a&\bigl(u^{I,h}, u^{II,h}\bigr) - f'(u^{I,h})\bigr) \, 
    (u^{II,h} - u^{I,h}) \, d V^{I,h}  
  - \bigl(a(u^I, u^{II}) - f'(u^I)\bigr) \, (u^{II} - u^I) \, d V^{I,h}  
  \\ 
& + 
  \bigl(a\bigl(u^I, u^{II}\bigr) - f'(u^I)\bigr) \, 
  (u^{II} - u^I) \, d V^{I,h}
- 
  \bigl(a\bigl(u_\eps^I, u_\eps^{II}\bigr) - f'(u_\eps^I)\bigr) \, 
  (u_\eps^{II} - u_\eps^I) \, d V^{I,h}\\ 
& + 
  \bigl(a\bigl(u_\eps^I, u_\eps^{II}\bigr) - f'(u_\eps^I)\bigr) \, 
  (u_\eps^{II} - u_\eps^I) \, d V^{I,h} 
- 
  \bigl(a\bigl(u_\eps^I, u_\eps^{II}\bigr) - f'(u_\eps^I)\bigr) \, 
  (u_\eps^{II} - u_\eps^I) \, d\bar V^I\\ 
& + 
  \bigl(a\bigl(u_\eps^I, u_\eps^{II}\bigr) - f'(u_\eps^I)\bigr) \, 
  (u_\eps^{II} - u_\eps^I) \, d\bar V^I
   - 
  \bigl(a\bigl(u^I, u^{II}\bigr) - f'(u^I)\bigr) \, 
  (u^{II} - u^I) \, d\bar V^I,  
\endaligned 
$$ 
which, with obvious notation, yields a decomposition for $T^h_{2,1}$  
$$
T_{2,1}^h = M_1^h+ M_2^h + M_3^h + M_4^h. 
\tag 3.39   
$$
Using (3.17) and the local convergence property (3.31), we obtain 
$$ 
\aligned 
|M_1^h| 
& \leq C \, \sum_{j=0}^n \int_{(y_{j}^h, y_{j+1})} 
  \bigl|d V^{I,h}\bigr| \,  
  \Big(\sup_{(y_{j}^h, y_{j+1})} |u^{I,h} - u^I| 
  + \sup_{(y_{j}^h, y_{j+1})} |u^{II,h} - u^{II}| \Big) \\ 
& \leq C \, \eps. 
\endaligned 
\tag 3.40  
$$ 
Similarly using (3.17) and (3.33) we obtain 
$$
\aligned 
|M_2^h| 
& \leq C \, \sum_{j=0}^n \int_{(y_{j}^h, y_{j+1})} 
  \bigl|d V^{I,h}\bigr| \,  
  \Big(\sup_{(y_{j}^h, y_{j+1})} |u^I - u^I_\eps| 
  + \sup_{(y_{j}^h, y_{j+1})} |u^{II} - u^{II}_\eps| \Big) \\ 
& \leq C \, \eps. 
\endaligned 
\tag 3.41
$$ 
Dealing with $M_4^h$ is similar:
$$
\aligned 
|M_4^h| 
& \leq C \, \sum_{j=0}^n \int_{(y_{j}^h, y_{j+1})} 
  \bigl|d \bar V^I\bigr| \,  
  \Big(\sup_{(y_{j}^h, y_{j+1})} |u^I - u^I_\eps| 
  + \sup_{(y_{j}^h, y_{j+1})} |u^{II} - u^{II}_\eps| \Big) \\ 
& \leq C \, \eps. 
\endaligned 
\tag 3.42  
$$ 

Finally to treat $M_3^h$ we observe that, since $u_\eps^I$ and $u_\eps^{II}$
are continuous functions on each interval $(y_{j}^h, y_{j+1})$ 
and since $d V^{I,h}$ is sequence 
of bounded measures converging weakly-star toward $d \bar V^I$, 
we have for all $h$ sufficiently small 
$$
|M_3^h| \leq \eps. 
\tag 3.43  
$$
Combining (3.39)--(3.43) we get 
$$
\bigl|T_{2,1}^h\bigr| \leq C \, \eps. 
\tag 3.44
$$
Combining (3.36)--(3.38) and (3.44) we obtain  
$$
\bigl|T_2^h\bigr| \leq C \, \eps
$$
and thus with (3.34)-(3.35) 
$$
\bigl|\Omega_2^h(\tau)\bigr| \leq C \, \eps \quad \text{ for all $h$ sufficiently small.} 
$$ 
Since $\eps$ is arbitrary, this completes the proof of (3.30).

%___________________________________________________________________________ 
\vskip.5cm 

{\bf Step 3 :} Consider now the term 
$$
\Omega_1(\tau) = \sum_{(x,\tau) \in \LL(a^h)\cap\JJ(u^{I,h})} 
 \Big(2 \, m + TV(b^h)\Big) \, 
 \bigl|a^h(x-,\tau) - \lam^h(x,\tau) \bigr| 
   \, |u^{II,h}(x-, \tau) - u^{I,h}(x-, \tau)|. 
\tag 3.45 
$$
On one hand, observe that 
$$
TV(b^h(\tau)) = TV(u^{I,h}(\tau)) + TV(u^{II,h}(\tau))  
\longrightarrow \bar V^I(+\infty, \tau) + \bar V^{II}(+\infty, \tau). 
\tag 3.46 
$$ 
For all but countably many $\tau$ the following holds. 
Extracting a subsequence if necessary we can always assume that 
for each $j$ either $(y_j^h, \tau) \in \LL(a^h)$ for all $h$, or 
else $(y_j^h, \tau) \notin \LL(a^h)$ for all $h$. 
Then consider the following three sets: denote by $J_1$ the set of 
indices $j$ such that $(y_j^h, \tau) \in \LL(a^h)$ and 
$(y_j, \tau) \in \LL(a)$. Let  $J_2$ the set of 
indices $j$ such that $(y_j^h, \tau) \notin \LL(a^h)$ and 
$(y_j, \tau) \in \LL(a)$. Finally $J_3$ is the set of 
indices $j$ such that $(y_j^h, \tau) \in \LL(a^h)$ and 
$(y_j, \tau) \notin \LL(a)$.

The local convergence property (3.27) implies  
$$
\aligned 
& \sum_{j \in J_1} 
\bigl|a\bigl(u^{I,h}(y^h_j-), u^{II,h}(y^h_j-)\bigr) 
 - a\bigl(u^{I,h}(y^h_j-), u^{I,h}(y^h_j+)\bigr) \bigr| 
   \, |u^{II,h}(y_j^h-) - u^{I,h}(y_j^h-)| \\ 
& \longrightarrow 
\sum_{j \in J_1} 
\bigl|a\bigl(u^I(y_j-), u^{II}(y_j-)\bigr) 
 - a\bigl(u^I(y_j-), u^I(y_j+)\bigr) \bigr| 
  \, |u^{II}(y_j-) - u^I(y_j-)|. 
\endaligned 
\tag 3.47
$$
(Indeed, given $\eps>0$, choose finitely many jump points as 
in (3.31) and use (3.27) with $\eps$ replaced with $\eps \, |u^I(z+) - u^I(z+)|$). 

On the other hand for indices in $J_2$ or $J_3$ we have 
$$ 
\sum_{j \in J_2 \cup J_3} 
\bigl|a\bigl(u^{I,h}(y^h_j-), u^{II,h}(y^h_j-)\bigr) 
 - a\bigl(u^{I,h}(y^h_j-), u^{I,h}(y^h_j+)\bigr) \bigr| 
   \, |u^{II,h}(y_j^h-) - u^{I,h}(y_j^h-)|
\longrightarrow
0
\tag 3.48
$$
but 
$$ 
\sum_{j \in J_2 \cup J_3} 
\bigl|a\bigl(u^I(y_j-), u^{II}(y_j-)\bigr) 
 - a\bigl(u^I(y_j-), u^I(y_j+)\bigr) \bigr| 
   \, |u^{II}(y_j-) - u^I(y_j-)|
=0. 
\tag 3.49
$$ 
Indeed, for each $j \in J_2$, 
$y_j$ is a Lax shock but $y_j^h$ is not.
Extracting a subsequence if necessary, it must be that the Lax inequalities are violated 
on the left or on the right side of $y_j^h$ for all $h$. 
So it must be that, assuming that it is the case on the left side,  
$a\bigl(u^I(y_j-), u^{II}(y_j-)\bigr) - a\bigl(u^I(y_j-), u^I(y_j+) \geq 0$ 
while 
$a\bigl(u^I(y_j^h-), u^{II}(y_j^h-)\bigr) - a\bigl(u^I(y_j^h-), u^I(y_j^h+) \leq 0$ 
for all $h$. But the latter converges toward the former by the local uniform 
convergence, which proves that 
$a\bigl(u^I(y_j-), u^{II}(y_j-)\bigr) - a\bigl(u^I(y_j-), u^I(y_j+) =0$.  

Combining (3.45)--(3.49) yields 
$$ 
\Omega_1^h \to 
\int_0^t \sum_{(x,\tau) \in \LL(a)} 
\bar q(\tau) \, \bigl|a(x-,\tau) - \lam(x,\tau) \bigr| \, 
|u^{II}(x, \tau) - u^I (x, \tau)| \, d\tau, 
\tag 3.50
$$ 
where $\bar q : = 2 \, m + \bar V^I(+\infty) + \bar V^{II}(+\infty)$. 

%___________________________________________________________________________ 
\vskip.5cm 

{\bf Step 4 : Continuity of the weighted norm.} 
  
Fix some time $t$.  
Recall that the weight $\bar w(t)$ is defined based on the total variation functions 
$\bar V^{II}$  and $\bar V^I$ and on the function $u^{II}(t) - u^I(t)$. 
The weight $w^h(t)$ is defined based on the total variation functions 
$V^{II,h}$ and $V^{I,h}$ and on the function $u^{II,h}(t) - u^{I,h}(t)$.
On the other hand, $u^{II,h} - u^{I,h}(t) \to u^{II} - u^I(t)$, 
$V^{II,h} \to\bar V^{II}$ and $V^{I,h} \to\bar V^I$.  
Therefore we have 
$$ 
w(x,t) = \bar w(x,t) \quad \text{ whenever } u^{II}(x,t) - u^I(x,t) \neq 0. 
\tag 3.51
$$

Combining (3.51) and the $L^1$ convergence $u^{II,h} - u^{I,h}(t) \to u^{II} - u^I(t)$, 
we have 
$$ 
\|u^{II}(t) - u^I(t)\|_{\bar w(t)} 
= 
\lim_{h \to 0} \|u^{II,h}(t) - u^{I,h}(t)\|_{w^h(t)}.  
\tag 3.52  
$$   

%___________________________________________________________________________ 
\vskip.5cm 

{\bf Step 5 : The left-hand side of $(3.16)$ is bounded above by $(3.21)$.} 

First of all, the inequality 
$$
\aligned 
& \int_0^t \int_\RR \bigl(a\bigl(u^I, u^{II}\bigr) - f'(u^I)\bigr) \, 
(u^{II} - u^I) \, dV^I(y, \tau) d\tau \\ 
& \leq 
\int_0^t \int_\RR \bigl(a\bigl(u^I, u^{II}\bigr) - f'(u^I)\bigr) \, 
(u^{II} - u^I) \, d\bar V^I(y, \tau) d\tau 
\endaligned  
\tag 3.53 
$$  
is a direct consequence of (3.19) and 
the definition of the nonconservative product in (3.15).

On the other hand, by the definition 
of the weighted norm and because of (3.20), 
similarly to (3.51) we have the inequality
$$ 
w(x,t) \leq \bar w(x,t) \quad \text{ whenever } u^{II}(x,t) - u^I(x,t) \neq 0. 
\tag 3.54
$$
Hence, by (3.53) and (3.54) 
the left-hand side of (3.16) is bounded above by (3.21).
This completes the proof of Theorem 3.5. $\quad\qed$   
\enddemo

%=========================================================================== 

\heading{4. Generalized Characteristics and Maximum Principle}
\endheading

We now return to the setting in Section 2 and aim at 
extending the analysis therein to arbitrary functions 
of bounded variation. For exact solutions of the hyperbolic equation
$$
\del_t \psi + \del_x \bigl( a \, \psi\bigr) = 0, 
\tag 4.1
$$ 
we will establish a maximum principle: Any solution of $(4.1)$ remains non-negative 
for all times if it is so initially. For a more precise (local) statement, our proof will 
make use of Dafermos-Filippov theory of generalized characteristics. 

Our main assumption throughout this section is the following: 
$$
\text{There exists a constant $E$ such that } 
\del_x a \leq {E \over t}. 
\tag 4.2 
$$
This is nothing but a generalization of the well-known 
Oleinik's entropy inequality. To motivate (4.2), let us recall the
following result. 

Let $f$ be a strictly convex function and $u$ be an
entropy solution (with bounded variation for all times) of 
the conservation law 
$$
\del_t u + \del_x f(u) = 0, \qquad u(x,t) \in \RR. 
\tag 4.3   
$$  
Then is is known that there exists a constant $C=C(u)$ such that 
$$
\del_x u \leq {C \over t}. 
\tag 4.4 
$$

\proclaim{Lemma 4.1} If $u^I$ and $u^{II}$ are two entropy solutions of 
the conservation law $(4.3)$, then the averaging speed
$$
a=a(u^I, u^{II}) := {f(u^{II}) - f(u^I) \over u^{II} - u^I}. 
\tag 4.5 
$$ 
satisfies our assumption $(4.2)$, with $E= \sup f''  \,  \bigl(C(u^I) + C(u^{II})\bigr) /2$. 
\endproclaim

\demo{Proof} Let us fix some time $t >0$. 
On each Borel set consisting of points of continuity 
of both $u^I$ and $u^{II}$, the following holds: 
$$
\aligned 
\del_x a 
& = \del_x \int_0^1 f'\bigl(\theta \, u^I + (1 - \theta) \, u^{II}\bigr) \, d\theta\\ 
& = \int_0^1 f''\bigl(\theta \, u^I + (1 - \theta) \, u^{II}\bigr) \, 
\bigl(\theta \, \del_x u^I + (1 - \theta) \, \del_x u^{II}\bigr) \, d\theta\\ 
& \leq
 \int_0^1 \sup f''  \, 
\bigl(\theta \, {C(u^I) \over t} + (1 - \theta) \, {C(u^{II}) \over t} \bigr) \, d\theta\\ 
& \leq \sup f''  \,  {C(u^I) + C(u^{II}) \over 2 \, t}.  
\endaligned 
$$ 

On the other hand, at a point $x$ where one of $u^I$ or $u^{II}$ is discontinuous, 
we have with an obvious notation 
$$
a_+ - a_- 
=  \int_0^1 f'\bigl(\theta \, u^I_+ + (1 - \theta) \, u^{II}_+\bigr) \, d\theta
- 
 \int_0^1 f'\bigl(\theta \, u^I_- + (1 - \theta) \, u^{II}_-\bigr) \, d\theta
 \leq 0, 
$$ 
since $f'$ is an increasing function and (for instance  by (4.4)) 
both $u^I$ and $u^{II}$ satisfy $u^I_+ \leq u^I_-$ and $u^{II}_+ \leq u^{II}_-$. 
$\quad\qed$  
\enddemo 

%==========================================================

By definition, a generalized characteristic $y= y(t)$ 
associated with the coefficient $a$ must satisfy for almost every $t$ 
(in its domain of definition) 
$$
a_+(y(t), t)   \leq    y'(t) \leq  a_-(y(t), t).  
\tag 4.6 
$$ 
According to Filippov's theory of differential equations \cite{\refFilippov}, 
through each point $(\bar x,\bar t)$ there pass a maximal and 
a minimal generalized characteristic.

\proclaim{Definition 4.2} A generalized characteristic is said to be 
genuine iff for almost every $t$ it satisfies 
$$
y'(t) \in \bigl\{a_-(y(t), t), a_+(y(t), t)\bigr\}. 
\tag 4.7
$$ 
\endproclaim

\proclaim{Proposition 4.3}  Any minimal backward generalized characteristic is  
genuine and for almost every $t$ satisfies 
$$
y'(t) = a_-(y(t), t). 
\tag 4.8
$$ 
Similarly, for a maximal backward generalized characteristic we have 
$y'(t) = a_+(y(t), t)$. 
\endproclaim

\demo{Proof} Here we only rely on the following consequence of (4.2): 
$a_+ \leq a_-$ at each discontinuity point of the function $a$. Geometrically, 
this condition prevents the existence of rarefaction-shocks in $a$. 
On the other hand, rarefaction centers (also prevented by (4.2) for $t>0$) 
could still be allowed for the present purpose. 

Consider $(\bar x,\bar t)\in (-\infty ,+\infty)\times (0,\infty)$, and let 
$y(t):=y(t;\bar x,\bar t)$ be the minimal backward characteristic through
$(\bar x,\bar t)$. We prove that it is genuine on its domain $(s,\bar t]$. 
We proceed as in \cite{\refDafermostwo}
 and assume by contradiction that there is a
measurable set $J,~\bar J\subset (s,\bar t]$ of positive Lebesgue measure,
and $\ve >0$ such that 
$$
a_-(y(t),t)-y'(t)>2\ve, \qquad t\in J. \tag 4.9 
$$
For each $t\in J$ there exists $\delta(t)>0$ with the property
$$
a_+(x,t)\geq a_-(y(t),t)-\ve, \qquad x\in (y(t)-\delta(t),~y(t)). \tag 4.10
$$
Finally, there is a subset $I\subset J$ with $\mu^* (I)>0$ (here $\mu^*$
denotes the outer measure) and $\bar\delta >0$ such that
$\delta(t)>\bar\delta$ for $t\in I$.

Let $\tau$ be a density point of $I$, with respect to $\mu^*$. Thus there
exists $\bar r,~0<\bar r<\bar t-\tau$, so that
$$
{{\mu^*(I\cap [\tau,\tau+r])}\over r}>{{2|\alpha |+\ve}\over{2|\alpha |+2\ve}}\,,
\qquad 0<r\leq\bar r, \tag 4.11
$$
where
$$
\alpha :=\inf\bigl\{a_+(x,t)-a_-(y(t),t):~s<t\leq\bar t,
~y(t)-\bar\delta\leq x<y(t)\bigr\}.
$$
Now take a point $y\in (y(\tau)-\bar\delta,~y(\tau))$ with the property
$y>y(\tau)-{1\over 2}\ve\bar r$, and consider a forward characteristic 
$z(\cdot)$ through $(y,\tau)$. We first observe that 
$$
z(t)<y(t), \qquad t>\tau,
$$
since $y(t)$ is the minimal backward characteristic through $(\bar x,\bar
t)$.

In addition, we have 
$$
z(t)>y(t)-\bar\delta, \qquad t>\in [\tau,~\tau+\bar r].
$$
Indeed, suppose by contradiction that for some $r\in (0,\bar r]$,
$z(t)>y(t)-\bar\delta$ for $t>\in [\tau,~\tau+r)$, but 
$z(\tau+r)=y(\tau+r)-\bar\delta$. Then
$$
\aligned
0&=z(\tau+r)-y(\tau+r)-\bar\delta=y+\int_\tau^{\tau+r} z'(t)dt~-~y(t)
~-~\int_\tau^{\tau+r} y'(t)dt~+~\bar\delta\\
&>\int_\tau^{\tau+r}(z'(t)-y'(t))dt\\
&=\int_{I\cap [\tau,\tau+r]}(z'(t)-a_-(y(t),t)+a_-(y(t),t)-y'(t))dt
\\
& \, 
+\int_{[\tau,\tau+r]\setminus I}(z'(t)-a_-(y(t),t)+a_-(y(t),t)-y'(t))dt\\
&\geq\ve\mu^*\bigl(I\cap [\tau,\tau+r]\bigr)+
 \alpha\bigr(r-\mu^*(I\cap[\tau,\tau+r])\bigr) >0,\\
\endaligned
$$
by (4.9)-(4.11), which leads to a contradiction. In the same way one obtains
$$
\aligned 
0&>z(\tau+\bar r)-y(\tau+\bar r)=y+\int_\tau^{\tau+\bar r} z'(t)dt~-~y(t)
~-~\int_\tau^{\tau+\bar r} y'(t)dt\\
&>\ve\mu^*\bigl(I\cap [\tau,\tau+\bar r]\bigr)+
 \alpha\bigr(\bar r-\mu^*(I\cap[\tau,\tau+r])\bigr)-{1\over 2}\ve\bar r >0,\\
\endaligned
$$
which gives another contradiction. For the maximal backward
characteristic the proof is similar.
$\quad\qed$  
\enddemo

\proclaim{Proposition 4.4} Forward characteristics leaving from some 
$(\bar x, \bar t)$ are unique when $\bar t >0$. 
\endproclaim 

\demo{Proof} Suppose there were two forward characteristics $y(\cdot)$ and
$z(\cdot)$ through $(\bar x, \bar t)$ with $y(\tau)<z(\tau)$ for some
$\tau>\bar t$. By (4.2) we have
$$
z'(\tau)-y'(\tau)\leq a_-(z(t),t)-a_+(y(t),t)
\leq C_{\bar t} \bigr(z(\tau)-y(\tau)\bigl). 
\tag 4.12
$$
Integrating (4.12) from $\bar t$ to $\tau$ one gets $z(\tau)-y(\tau)=0$, which gives a
contradiction.
$\quad\qed$  
\enddemo

\proclaim{Theorem 4.5} Let $\psi=\psi(x,t)$ be a solution of $(4.1)$ 
such that on some interval $[\xi_0,\zeta_0]$ we have 
$$
\psi(x,0) \geq 0, \qquad x \in [\xi_0,\zeta_0].
\tag 4.13
$$  
Let $\xi=\xi(t)$ be any forward generalized characteristic leaving 
from $(\xi_0,0)$, and 
$\zeta=\zeta(t)$ be any forward generalized characteristic leaving 
from $(\zeta_0,0)$. 

Then we have for all $t \geq 0$ 
$$
\psi(x,t) \geq 0, \qquad x \in \bigl(\xi(t), \zeta(t)\bigl). 
\tag 4.14 
$$   
\endproclaim

Note that it may happen that $\xi(t)=\zeta(t)$ for $t$ large enough.

\demo{Proof} Observe that the two characteristics cannot cross and fix
any time $t>0$ such that $\xi(t) < \zeta(t)$. 
Fix also any two points such that $\xi(t) < \bar y < \bar z < \zeta(t)$. 
Let $y(t)$ and $z(t)$ be the maximal and minimal backward characteristics
emanating from $\bar y$ and $\bar z$, respectively. These characteristics 
can not leave the region limited by $\xi(t)$ and $\zeta(t)$.

Integrating (4.1) in the domain bounded by the characteristics $y(t)$ and $z(t)$, 
and using that these characteristics are genuine, so that the flux 
terms along the vertical boundaries vanish identically, we arrive at 
$$
\int_{\bar y}^{\bar z} \psi(x,t) \, dx 
= 
\int_{y(0)}^{z(0)} \psi(x,0) \, dx 
\geq 0. 
\tag 4.15
$$ 
The last inequality is due to the fact that $\psi(.,0) \geq 0$
and the inequalities $\xi_0 =\xi(0) \leq y(0)  \leq z(0) \leq \zeta(0)=\zeta_0$.
Since $\bar y$ and $\bar z$ are arbitrary, we obtain (4.14). 
$\quad\qed$  
\enddemo

%=========================================================================== 

\heading{5. A Sharp $L^1$ Estimate for Hyperbolic Linear Equations}
\endheading 

Based on the maximum principle established in Section 4,  
we now derive a sharp estimate for the weighted norm introduced 
in Section 2. 
We restrict attention again to the situation where 
$u^I$ and $u^{II}$ are two entropy solutions of 
the conservation law $(4.3)$ and $a$ is the averaging speed
given in (4.5). We define a weight by analogy with what was done 
in Section 2 in the special case of piecewise constant solutions. 

Given a solution $\psi$ of the equation (4.1), we introduce 
weighted $L^1$ norm in the following way. Set 
$$
V^I(x,t) = TV_{-\infty}^x (u^I(t)),   \qquad 
V^{II}(x,t) = TV_{-\infty}^x (u^{II}(t))  
\tag 5.1
$$ 
and fix some parameter $m \geq 0$. Then 
consider the weight-function defined, for each $t \geq 0$ 
and each point of continuity $x$ for $u^I(t)$ and $u^{II}(t)$, by 
$$
w(x,t) = \cases
m + V^I(\infty, t) - V^I(x,t)  
         + V^{II}(x,t) & \text {if } \psi(x,t) >0, \\
\\ 
m + V^I(x,t)
         + V^{II}(\infty, t) - V^{II}(x,t)   
         & \text {if } \psi(x,t) \leq 0.\\
\endcases 
\tag 5.2   
$$ 
It is immediate to see that 
$$ 
m \leq w(x,t) \leq m + TV(u^I(t)) + TV(u^{II}(t)), \qquad x \in \RR. 
\tag 5.3  
$$ 
Finally the weighted norm on the solutions $\psi$ of (4.1) is defined by 
$$
\|\psi(t) \|_{w(t)} := \int_\RR |\psi(x,t)| \, w(x,t) \, dx. 
$$ 
Note that the weight depends on the fixed solutions $u^I$ and $u^{II}$, 
but also on the solution $\psi$.

Our sharp estimate will involve the nonconservative product
$$  
\mu^I_\psi(t) = \bigl(a - f'(u^I(t))\bigr) \, \psi(t) \, dV^I(t) 
$$
defined for all almost every $t \geq 0$ by 
\roster
\item If $B$ is a Borel set included in the set of continuity points of $u^I(t)$
then 
$$
\mu_\psi^I(t) (B) = \int_B  \bigl(a(t) - f'(u^I(t)) \bigr) \, \psi(t) \, dV^I(t), 
\tag 5.4a 
$$
where the integral is defined in a classical sense; 
\item If $x$ is a point of jump of $u^I(t)$, then 
$$
\mu_\psi^I(t)(\{x\}) 
= \bigl(a(x-, t) - \lambda^I(x, t) \bigr) \, \psi( x-, t) 
\, \bigl|u^I(x+, t) - u^I(x-, t)\bigr|.  
\tag 5.4b
$$ 
\endroster 

Here $\lambda^I(x, t)$ is a the shock speed of the discontinuity in $u^I$ 
located at $(x,t)$. The measure $\mu_\psi^{II}(t)$ is defined similarly. 
Regarding the expression (5.4b), it is worth noting
that if $(x,t)$ is a point of approximate jump of $u^I$ and $\psi$, 
then the jump relation for the equation (4.1) reads 
$$
\bigl(a(x-, t) - \lambda^I(x, t) \bigr) \, \psi( x-, t)
= 
\bigl(a(x+, t) - \lambda^I(x, t) \bigr) \, \psi( x+, t). 
\tag 5.5
$$
In the same way we define
$$  
\mu^{II}_\psi(t) = \bigl(f'(u^{II}(t)) -a\bigr) \, \psi(t) \, dV^{II}(t). 
$$ 
We now prove: 

%__________________________________________________________________ 

\proclaim{Theorem 5.1} Let $u^I$ and $u^{II}$ be two entropy solutions 
of $(1.1)$ such that $u^{II}-u^I$ admits finitely many changes of sign. 
Let $\psi$ be any solution of bounded variation of the 
hyperbolic equation $(4.1)$ satisfying the constrain 
$$
\psi \,  \bigl( u^{II} - u^{I}\bigr) \geq 0. 
\tag 5.6 
$$
Then for all $0 \leq s \leq t$   
$$ 
\aligned 
\|\psi(t)\|_{w(t)} 
&+\int_s^t \sum_{(x,\tau) \in \LL(a)} 
\Big(2 \, m + TV(a) \Big) \, 
\bigl|a_-(x,\tau) - \lam(x,\tau) \bigr| \, 
|\psi_-(x, \tau)| \, d\tau \\ 
& +
\int_s^t \int_\RR \bigl(a(\tau) - f'(u^I(\tau)) \bigr) \, \psi(\tau) \, 
dV^I(\tau) d\tau 
+
\int_s^t \int_\RR \bigl(a(\tau) - f'(u^{II}(\tau)) \bigr) \, \psi(\tau) \, 
dV^{II}(\tau)  d\tau \\ 
& \leq 
 \|\psi(s)\|_{w(s)}. 
\endaligned 
\tag 5.7 
$$ 
\endproclaim

The assumption (5.6) is clearly satisfied with the choice 
$\psi =  u^{II} - u^{I}$. Therefore our previous result in Theorem 3.5 
(derived via a completely different proof) can be regarded as a corollary of 
Theorem 5.1.  

It is interesting to observe that, when $u^{II}=u^I$, the weight (5.2)
becomes constant, and therefore (5.7) reduces to the $L^1$ estimate. 
$$
\aligned 
\|\psi(t)\|_{L^1(\RR)} 
& + 
\int_s^t \sum_{(x,\tau) \in \LL(a)} 
\Big(2 \, m + TV(a) \Big) \, \bigl|a_-(x,\tau) - \lam(x,\tau) \bigr| \, 
|\psi_-(x, \tau)| \, d\tau \\ 
& \leq 
 \|\psi(s)\|_{L^1(\RR)}. 
\endaligned  
$$ 
Also, note that under the assumption (5.6) $\mu_\psi^I(t)$ and $\mu_\psi^{II}(t)$  
are positive except at points 
$(x,t)\in \LL(a)\cup\RS(a)$.
However, these negative terms are offset in (5.7) by the positve terms under the first 
integral.

%__________________________________________________________________________ 

\demo{Proof}  
Fix any positive time $t$. By assumption we have finitely many points 
$-\infty=y_0<y_1<\ldots<y_n<y_{n+1}=+\infty$ such
that, on each interval $(y_i,y_{i+1})$, we have $\psi (t)\geq 0$ when $i$ is odd and
$\psi (t)\geq 0$ when $i$ is even. For every $i=1, \cdots, n$,
consider the (unique by Proposition 4.4) forward characteristic 
$y_i(\cdot)$ associated with the coefficient $a$ 
and issuing from the initial point $(y_i,t)$.

We will focus attention on some interval $(y_i,y_{i+1})$ with $i$ odd, say, 
and with $-\infty<y_i<y_{i+1}<+\infty$. 
Except when specified differently, all of the characteristics to be 
considered from now on are associated with the solution $u^{II}$. 
For definiteness we will first study the case that
the forward characteristic $\chi_0(\cdot)$ (associated with $u^{II}$ and) 
issuing from the point $(y_i,t)$ is located on the right-side of the 
curve $y_i$, that is,  
$$
y_i(\tau) \leq \chi_0(\tau), \qquad t \leq \tau \leq t +\delta 
$$
for some $\delta>0$ sufficiently small. 

Fix some (sufficiently small) $\eps>0$ and denote by 
$y_i<z_1< \ldots <z_N<y_{i+1}$ the points where $u^I$ has 
a jump larger or equal to $\eps$, that is,  
$$
u_-^{II}(z_I,t)- u_+^{II}(z_I,t) \geq \eps, \qquad I=1,\ldots,N.
\tag 5.8 
$$ 
For each $I=1,\ldots,N$, 
consider also the forward characteristic $\chi_I(\cdot)$ issuing from the
point $(z_I,t)$. 
For definiteness, we will also assume 
that the forward characteristic $\chi_{N+1}(\cdot)$
issuing from $(y_{i+1},t)$ satisfies 
$$
\chi_{N+1}(\tau) \leq y_{i+1}(\tau), \qquad t \leq \tau \leq t + \delta
$$ 
for some $\delta>0$ sufficiently small.

Next, let us select a time $s>t$ with $s-t$ so small that the following 
properties hold:
\roster
\item"(a)" No intersection among the characteristics
$y_i,\chi_0,\chi_1,\ldots,\chi_N,\chi_{N+1},y_{i+1}$ occurs in the time
interval $[t,s]$.

\item"(b)" For $I=1,\ldots,N$, let $\zeta_I(\cdot)$ and $\xi_I(\cdot)$ be the
minimal and the maximal backward characteristics emanating from the point 
$(\chi_I(s),s)$. Then the total variation of $u^{II}(\cdot,t)$ over the
intervals $(\zeta_I(t),z_I)$ and $(z_I,\xi_I(t))$ should not exceed $\eps \over N$.

\item"(c)" Let $\zeta_0(\cdot)$ be the minimal backward characteristic 
emanating from $(y_i(s),s)$ and $\xi_0(\cdot)$ be the maximal backward 
characteristic emanating from $(\chi_0(s),s)$. Then the total variation of
$u^{II}(\cdot,t)$ over the intervals $(y_i,\xi_0(t))$ and $(\zeta_0(t),y_i)$
should not exceed $\eps$.

\item"(d)"  Let $\zeta_{N+1}(\cdot)$ be the minimal backward characteristic 
emanating from the point $(\chi_{N+1}(s),s)$ and $\xi_{N+1}(\cdot)$ be the maximal backward 
characteristic emanating from $(y_{i+1}(s),s)$. Then the total variation of 
$u^{II}(\cdot,t)$ over the intervals $(\zeta_{N+1}(t),y_{i+1})$ and
$(y_{i+1},\xi_{N+1}(t))$ should not exceed $\eps$.

\endroster 

For $I=0,\ldots,N$, and some integer $k$ to be fixed later, 
consider a mesh of the form  
$$
\chi_I(s) = x_I^0<x_I^1<\ldots<x_I^k<x_I^{k+1}=\chi_{I+1}(s).
\tag 5.9 
$$
For $I=0,\ldots,N$ and $j=1,\ldots,k$, consider also the maximal backward
characteristic $\xi_I^j(\cdot)$ emanating from the point $(x_I^j,s)$ and 
identify its intercept $z_I^j=\xi_I^j(t)$ by the horizontal line at time $t$. 
Finally set also 
$$
z_0^0 = y_i, \quad 
z_N^{k+1} = y_{i+1}, 
\quad 
z_{I-1}^{k+1} = z_I^0 = z_I, \qquad I=1,\ldots,N.
$$

%_____________________________________________________________________________  

To start the proof, we integrate the equation (4.1) satisfied by the function $\psi$, 
successively in each domain limited by the characteristics introduced above. 
Applying Green's theorem, we arrive at the following five formulas:

$(i)$ Integrating (4.1) on the region 
$$
\bigl\{(x,\tau) \, / \, t<\tau<s, \quad y_i(\tau)<x<\chi_0(\tau)\bigr\}
$$
and multiplying by $V^{II}(y_i,t)$ one gets 
$$
\aligned
\int_{y_i(s)}^{\chi_0(s)} \psi(x,s) \, V^{II}(y_i,t) \, dx
& +\int_t^s (y_i'-a_+) \, \psi_+(y_i(\tau),\tau) \, V^{II}(y_i,t) \, d\tau\\ 
& +\int_t^s (a_- -\lambda_0) \, \psi_-(\chi_0(\tau),\tau) \, V^{II}(y_i,t) \, d\tau =0.\\
\endaligned 
\tag 5.10i
$$

$(ii)$ Integrating (4.1) on each of the regions  
$$
\bigl\{(x,\tau) \, / \, t<\tau<s, \quad \xi_I^j(\tau)<x<\xi_I^{j+1}(\tau)\bigr\} 
$$
for $I=0,\ldots,N$ and $j=1,\ldots,k$, and then multiplying by $V^{II}(z_I^j+,t)$, 
one gets
$$
\aligned
&    \int_{x_I^j}^{x_I^{j+1}}   \psi(x,s) \, V^{II}(z_I^j+,t) \, dx
   - \int_{z_I^j}^{z_I^{j+1}}   \psi(x,t) \, V^{II}(z_I^j+,t) \, dx\\
&  + \int_t^s (\lambda_I^j-a_+) \psi_+(\xi_I^j(\tau), \tau) \, V^{II}(z_I^j+,t) \, d\tau
   + \int_t^s (a_- -\lambda_I^{j+1}) \, \psi_-(\xi_I^{j+1}(\tau), \tau) \, 
          V^{II}(z_I^j+,t) \, d\tau = 0. 
\endaligned 
\tag 5.10ii 
$$

$(iii)$ Integrating (4.1) on each of the regions  
$$
\bigl\{(x,\tau) \, / \, t<\tau<s, \quad \chi_I(\tau)<x<\xi_I^{1}(\tau)\bigr\} 
$$
for $I=0,\ldots,N$, and multiplying by $V^{II}(z_I+,t)$ one gets
$$
\aligned
&   \int_{\chi_I(s)}^{x_I^{1}} \psi(x,s) \, V^{II}(z_I+,t) \, dx
  - \int_{z_I}^{z_I^{1}} \, \psi(x,t) \, V^{II}(z_I+,t) \, dx\\
& + \int_t^s (\lambda_I-a_+) \, \psi_+(\chi_I(\tau), \tau) \, V^{II}(z_I+,t) \, d\tau
  + \int_t^s (a_- -\lambda_I^1) \, \psi_-(\xi_I^{1}(\tau), \tau) \, 
    V^{II}(z_I+,t) \, d\tau = 0. \\
\endaligned 
\tag 5.10iii
$$

$(iv)$ Integrating (4.1) on the regions 
$$
\bigl\{(x,\tau) \, / \, t<\tau<s, \quad \xi_I^k(\tau)<x<\chi_{I+1}(\tau)\bigr\} 
$$
for $I=0,\ldots,N$, and multiplying by $V^{II}(z_I^k+,t)$ one gets
$$
\aligned
&  \int_{x_I^k}^{\chi_{I+1}(s)} \, \psi(x,s) \, V^{II}(z_I^k+,t) \, dx
  -\int_{z_I^k}^{z_{I+1}} \psi(x,t) \, V^{II}(z_I^k+,t) \, dx\\
& +\int_t^s (\lambda_I^k-a_+) \, \psi_+(\xi_I^k(\tau), \tau) \, V^{II}(z_I^k+,t) \, d\tau
  +\int_t^s (a_- -\lambda_{I+1}) \, \psi_-(\chi_{I+1}(\tau), \tau) \, V^{II}(z_I^k+,t) \, d\tau = 0.\\
\endaligned 
\tag 5.10iv
$$

$(v)$ Finally integrating (4.1) on the last region 
$$
\bigl\{(x,\tau) \, / \, t<\tau<s, \quad \chi_{N+1}(\tau)<x<y_{i+1}(\tau)\bigr\} 
$$
and multiplying by $V^{II}(y_{i+1},t)$ one gets
$$
\aligned
&  \int_{\chi_{N+1}(s)}^{y_{i+1}(s)} \psi(x,s) \, V^{II}(y_{i+1},t) \, dx
+\int_t^s (\lambda_{N+1}-a_+) \, \psi_+(\chi_{N+1}(\tau), \tau) \, V^{II}(y_{i+1},t) \, d\tau \\
& +\int_t^s (a_- -y_{i+1}') \, \psi_-(y_{i+1}(\tau), \tau) \, V^{II}(y_{i+1},t) \, d\tau = 0. \\
\endaligned 
\tag 5.10v 
$$

Next, summing all of the formulas (5.10) leads us to the general identity: 
$$
\aligned
& \int_{y_i(s)}^{\chi_0(s)}  \psi(x,s) \, V^{II}(y_i,t) \, dx
  + \sum_{I=0}^N \sum_{j=0}^k\int_{x_I^j}^{x_I^{j+1}} \psi(x,s) \, V^{II}(z_I^j+,t) \, dx\\
& +\int_{\chi_{N+1}(s)}^{y_{i+1}(s)} \psi(x,s) \, V^I(y_{i+1},t) \, dx
 -\sum_{I=0}^N \sum_{j=0}^k \int_{z_I^j}^{z_I^{j+1}} \psi(x,t) \, V^{II}(z_I^j+,t) \, dx\\
& = 
-\sum_{I=0}^N \sum_{j=1}^k
\int_t^s [V^{II}(z_I^j+,t)-V^{II}(z_I^{j-1}+,t)] \, 
(\lambda_I^j-a_-) \, \psi_-(\xi_I^j(\tau),\tau) \, d\tau \\
&-\sum_{I=0}^N\int_t^s [V^{II}(z_I+,t)-V^{II}(z_{I-1}^{k}+,t)] \, 
(\lambda_I^j-a_-) \, \psi_-(\chi_I(\tau),\tau) \, d\tau \\
&-\int_t^s[V^{II}(y_i+,t)-V^{II}(y_i,t)] \, 
(\lambda_0-a_-) \, \psi_-(\chi_0(\tau),\tau) \, d\tau \\
&-\int_t^s[V^{II}(y_{i+1}+,t)-V^{II}(z_N^k,t)] \, 
(\lambda_{N+1}-a_-) \, \psi_-(\chi_{N+1}(\tau),\tau) \, d\tau \\
&-\int_t^s (y_i'-a_+) \, \psi_+(y_i(\tau),\tau) \, V^{II}(y_i,t) \, d\tau 
-\int_t^s (a_- -y_{i+1}') \, \psi_-(y_{i+1}(\tau),\tau) \, V^{II}(y_{i+1},t) \, d\tau. \\
\endaligned 
\tag 5.11  
$$
To estimate the right-hand side of (5.11), we recall that 
the solution $u^I$ of a scalar conservation laws satisfies
$$
V^{II}(y_i,t)\geq V^{II}(\chi_0(s),s), \qquad 
V^{II}(z_I^j+,t)\geq V^{II}(x_I^j+,s),
$$
for $I=0,\ldots,N$ and $j=0,\ldots,k$. Hence, choosing the difference 
$x_I^{j+1}-x_I^j$ in (5.9) sufficiently small and since the function 
$V^{II}(\cdot,t)$ is nondecreasing, we conclude that the left-hand side of (5.11) 
can be bounded from below, as follows: 
$$
\text{ L.H.S. } 
\geq \int_{y_i(s)}^{y_{i+1}(s)}\psi(x,s)V^{II}(x,s)dx - (s-t)\ve
  -\int_{y_i(t)}^{y_{i+1}(t)} \psi(x,t)V^{II}(x,t)dx.
\tag 5.12
$$

Estimating the right-hand side of (5.11) is more involved. First note 
that each term arising in the left-hand 
side of (5.11) is non-positive. This follows from our condition (5.6). 
Indeed, 
consider a point $(x,s)$ of approximate jump or approximate continuity 
of $u^I$, $u^{II}$ and $\psi$. If all of these functions are continuous, 
the result is trivial. Call $\lambda$ the discontinuity speed. 
Based on the jump relation (5.5), 
we see that either $\psi_- \, (\lambda - a_-) = \psi_+ \, (\lambda - a_+) = 0$, 
or else 
all of the terms $\psi_-$, $\lambda - a_-$, $\psi_+$, and $\lambda - a_+$
are distinct from zero. 

Suppose first that $(x,s)$ is a point in the interior of the region limited by 
the two curves $y_i(.)$ and $y_{i+1}(.)$. 
In the latter case, since $\psi \geq 0$ in the region under consideration, 
we deduce that $\psi_- >0$ and $\psi_+ > 0$, while  
the terms $\lambda - a_-$ and $\lambda - a_+$ are 
either both negative or both positive. Actually, 
in view of the sign condition (5.6), we have $u^{II}_\pm - u^I_\pm \geq 0$ 
and, therefore, $\lambda - a_\pm \geq$ as follows from (3.10) (here we are dealing 
with a jump of $u^{II}$). 

Consider next a point of the boundary $y_i$, for instance. 
So we now have  $\psi_- <0$ and $\psi_+ > 0$, while  
the terms $\lambda - a_-$ and $\lambda - a_+$ opposite sign. 
Since no rarefaction-shock can arise, the discontinuity must be a Lax shock 
and so $\lambda - a_-<0$ and $\lambda - a_+>0$. Again the corresponding 
term in (5.11) has a favorable sign.  
(Observe that the condition (5.6) was not used in this second case.)

%___________________________________________________________________ 

Then, for all $I=0,\ldots,N$ and $j=1,\ldots,k$, let 
$\theta_I^j(\cdot)$ be the (maximal, for definiteness) backward characteristic 
associated with $u^{I}$ and issuing from the point $(\xi_I^j(\tau),\tau)$. 
Denote also by $\theta(z_I^j;\tau)$ its intercept with the horizontal line
at time $t$. Setting 
$$
\tilde a(x,t;\tau) := {{f(u^{II}(x,t))-f(u^{I}(\theta(x,\tau),t))}
\over {u^{II}(x,t)-u^{I}(\theta(x;\tau),t)}}
$$
and using that the solution $u^{I}$ remains constant along the 
characteristic $\theta_I^j(\cdot)$,
we obtain 
$$
(\lambda_I^j-a_-)(\xi_I^j(\tau)) = \lambda_I^j(z_I^j) - \tilde a(z_I^j,t;\tau).
\tag 5.13
$$
Then consider the (maximum, for definiteness) backward characteristic 
$y_I^j(\cdot)$ 
associated with $a$ and issuing from the point $(\xi_I^j(\tau),\tau)$. 
By integrating $\psi$ along the characteristic $y_I^j(\cdot)$ and using
the inequality (4.4), we arrive at a lower bound for $\psi$ 
$$
\psi(\xi_I^j(\tau),\tau)\geq\psi(y_I^j(t),t)\left( {t\over\tau}\right)^E, 
\qquad t < \tau < s.
\tag 5.14 
$$
Upon choosing $x_I^{j+1}-x_I^j$ in (5.9) so small that the oscillation of 
$V^{II}_c(\cdot)$ over each interval $(z_I^j-z_I^{j+1})$ does not exceed
$\ve$ and recalling the standard estimates on Stieltjes integrals we 
deduce from (5.11)-(5.13) that 
$$ 
\aligned
&\sum_{I=0}^N \sum_{j=1}^k
\int_t^s[V^{II}(z_I^j+,t)-V^{II}(z_I^{j-1}+,t)] \, 
\, \bigl((\lambda_I^j-a_-) \, \psi_-(\xi_I^j(\tau),\tau) \bigr) \, d\tau \\
&\geq\sum_{I=0}^N \sum_{j=1}^k
\int_t^s[V^{II}_c(z_I^j,t)-V^{II}_c(z_I^{j-1},t)] \, 
\bigl(\lambda_I^j(z_I^j)-\tilde a(z_I^j,t,\tau)\bigr)
\psi(y_I^j(t),t) \, \left( {t\over\tau}\right)^E \, d\tau \\
&
\geq\int_t^s\sum_{I=0}^N \sum_{j=1}^k \left( \int_{z_I^{j-1}}^{z_I^j}
\bigl(\lambda_I^j(x)-\tilde a(x,t,\tau)\bigr) \, 
\psi(x,t) \,  dV^{II}_c(x,t) - c \, \ve\right) 
\left( {t\over\tau}\right)^E   d\tau \\
&
=\int_t^s \left( \int_{y_i}^{y_{i+1}} 
\bigl(\lambda_I^j(x)-\tilde a(x,t,\tau)\bigr)
\psi(x,t) \,  dV^{II}_c(x,t) - c(y_{i+1}-y_i)\ve\right)
\left( {t\over\tau}\right)^E   d\tau. \\
\endaligned 
\tag 5.15 
$$

%___________________________________________________________________ 

We now combine (5.10), (5.11) and (5.15), divide the resulting
inequality by $s-t$, and let $s\searrow t$, $\eps \to 0$, obtaining
the following inequality: 
$$
\aligned
{{d^+ \over dt}}\int_{y_i(t)}^{y_{i+1}(t)} \psi(x,t)V^{II}(x,t)dx 
  \leq
& -\int_{y_i}^{y_{i+1}} \bigl(\lambda_I^j(x)-a(x,t)\bigr) \, \psi(x,t)
  \, dV^{II}_c(x,t)\\
& -\sum_{(x,t) \in \JJ(u^{II})} 
  \bigl(u^{II}_-(x,t) - u^{II}_+(x,t)\bigr) \, (\lam^I -a_-)(x,t) \, \psi_-(x,t)  \\
& -\bigl(u^{II}_-(y_i,t) - u^{II}_+(y_i,t)\bigr)\,(\lam^I -a_-)(y_i,t) \,
\psi_-(y_i,t) 
 \\
& -\bigl(u^{II}_-(y_{i+1},t) - u^{II}_+(y_{i+1},t)\bigr)\,
    (\lam^I -a_-)(y_{i+1},t)\, \psi_-(y_{i+1},t)    \\
& - (y_i'-a_+) \, \psi_+(y_i,t)  \, V^{II}(y_i,t) \\
& - (a_- -y_{i+1}') \, \psi_-(y_{i+1},t)  \, V^{II}(y_{i+1},t) \\
\endaligned 
\tag 5.16 
$$
The third and fourth terms in the right-hand side of (5.16) 
are due to the fact that $\chi_0$ and $\chi_{N+1}$
lie inside the region limited by $y_i$ and $y_{i+1}$.

We can next focus on the intervals $(y_i,y_{i+1})$ with $i$ even. 
Based on a completely symmetric argument and using now 
the weight $m+V^{II}(\infty,t)-V^{II}(\cdot,t)$ instead of 
$V^{II}(\cdot,t)$, we obtain
$$
\aligned
{{d^+}\over{dt}}&\int_{y_i(t)}^{y_{i+1}(t)} 
\bigl(-\psi(x,t)\bigr)\,\Big(m+V^{II}(\infty,t)-V^{II}(x,t)\Big)dx  \\
\leq & \int_{y_i}^{y_{i+1}}
\bigl(\lambda_I^j(x)-a(x,t)\bigr)
\bigl(-\psi(x,t)\bigr) dV^{II}_c(x,t) \\
&+\sum_{(x,t) \in \JJ(u^{II})} 
\bigl(u^{II}_-(x,t) - u^{II}_+(x,t)\bigr)\,(\lam^I -a_-)(-\psi_-)(x,t)  \\
&+\bigl(u^{II}_-(y_i,t) - u^{II}_+(y_i,t)\bigr)\,(\lam^I -a_-)(-\psi_-)(y_i,t) \\
&+\bigl(u^{II}_-(y_{i+1},t) - u^{II}_+(y_{i+1},t)\bigr)\, (\lam^I -a_-) 
  (-\psi_-)(y_{i+1},t) \\
&-(y_i'-a_+)(-\psi_+)(y_i,t)\Big(m+V^{II}(\infty,t)- V^{II}(y_i,t)\Big) \\
&-(a_- -y_{i+1}')(-\psi_-)(y_{i+1},t)\Big(m+V^{II}(\infty,t)-V^{II}(y_{i+1},t)\Big). \\
\endaligned 
\tag 5.17
$$
By summation over $i=1,\ldots,n$ in (5.16) for $i$ odd and in 
(5.17) for $i$ even respectively, we obtain
$$
\aligned
{{d^+}\over{dt}}&\int_{-\infty}^{+\infty} 
\bigl[\psi(x,t)\bigr]^+V^{II}(x,t)+\bigl[-\psi(x,t)\bigr]^+
\Big(m+V^{II}(\infty,t)-V^{II}(x,t)\Big)dx \\
\leq&-\sum_{(x,t) \in \LL(a)\cap\JJ(u^{II})} 
\Big(m + V^{II}(\infty,t) \Big)\,\bigl|\lam(x,t)-a_-(x,t)\bigr|\,|\psi_-(x,t)|\\
&-\sum_{(x,t) \in \JJ(u^{II})} 
\bigl(u^{II}_-(x,t) - u^{II}_+(x,t)\bigr)\,\bigl(\lam^I(x,t)-a_-(x,t)\bigr)\,
\psi_-(x,t) \\
&-\int_\RR \bigl(f'(u^{II}(y,t)) - a(y,t)\bigr)\,\psi(y,t)\,dV_c^{II}(y,t), \\
\endaligned 
\tag 5.18
$$ 
where the superscript $+$ denotes the positive part of the functions $\psi$
and $-\psi$ respectively.

Consider now the case where
$$
\chi_0(\tau)\leq y_i(\tau) , \qquad t \leq \tau \leq t +\delta ,
$$
and
$$
y_{i+1}(\tau)\leq \chi_{N+1}(\tau), \qquad t \leq \tau \leq t + \delta .
$$
Assume that there exists a time $\bar\tau>t$ such that 
$$
\chi_0(\bar\tau)<y_i(\bar\tau), \qquad\qquad y_{i+1}(\bar\tau)<\chi_{N+1}(\bar\tau)
$$
(otherwise the curves of the two pairs will coincide, and we can reduce to the 
previous 
case). Let now $\xi_0(\cdot)$ be the maximal backward characteristic emanating 
from $(y_i(\bar\tau),\bar\tau)$, and $\zeta_{N+1}(\cdot)$ be the minimal backward 
characteristic emanating from the point $(y_{i+1}(\bar\tau),\bar\tau)$. 
Since characteristics cannot cross, we have that
$$
y_i(t)<\xi_0(t),~\zeta_{N+1}(t)<y_{i+1}(t),
$$
Then, by finite propagation speed, there exists a time $s>t$ such that
$$
y_i(\tau)<\xi_0(\tau),~\zeta_{N+1}(\tau)<y_{i+1}(\tau),\qquad t\leq\tau<s,
$$
$$
y_i(s)=\xi_0(s),~\zeta_{N+1}(s)=y_{i+1}(s).
$$
Instead of properties (c), (d), we will require that $s$ satisfies the
following:

\roster

\item"(c')" Let $\zeta_0(\cdot)$ be the minimal backward characteristic 
emanating from $(\chi_0(s),s)$. Then the total variation of
$u^{II}(\cdot,t)$ over the intervals $(y_i,\xi_0(t))$ and $(\zeta_0(t),y_i)$
should not exceed $\eps$. 

\item"(d')"  Let $\xi_{N+1}(\cdot)$ be 
the maximal backward characteristic emanating from $(\chi_{N+1}(s),s)$. 
Then the total variation of $u^{II}(\cdot,t)$ over the intervals 
$(\zeta_{N+1}(t),y_{i+1})$ and $(y_{i+1},\xi_{N+1}(t))$ should not 
exceed $\eps$. 

\endroster 
{}From then on we can proceed as before.
Finally we write the inequality in (5.18) exchanging the roles of $u^I$ and 
$u^{II}$, and combining it with (5.17) we arrive exactly at the desired 
inequality (5.7) and the proof of Theorem 5.1 is completed.
$\quad\qed$  
\enddemo

%=========================================================================== 
\heading{Acknowledgements} 
\endheading

The authors are very grateful to C. Dafermos who communicated to them 
his lecture notes on the Liu-Yang functional in the context of 
general functions with bounded variation.

%=========================================================================== 
\heading{References} 
\endheading

\item{[\refBressan]} Bressan A., 
{\it Hyperbolic Systems of Conservation Laws,\/} Oxford Univ. Press, 
to appear. 

\item{[\refBCP]}     Bressan A., Crasta G. and Piccoli B., 
Well-posedness of the Cauchy problem for $n\times n$ systems 
of conservation laws, Mem. Amer. Math. Soc., to appear. 

\item{[\refBL]}      Bressan A. and LeFloch P.G., 
Structural stability and regularity of entropy solutions to 
systems of conservation laws, Indiana Univ. Math. J. 48 (1999), 43--84.   

\item{[\refBLY]}     Bressan A., Liu T.P. and Yang T., 
$L^1$ stability estimate for $n \times n$ conservation laws, 
Arch. Rational Mech. Anal. 149 (1999), 1--22.

\item{[\refCLone]}      Crasta G. and LeFloch P.G., 
Existence theory for a class of strictly hyperbolic systems, 
in preparation. 

\item{[\refCLtwo]}      Crasta G. and LeFloch P.G.,  in preparation. 

\item{[\refDafermosone]} Dafermos C.M., 
Polygonal approximations of solutions of the initial 
value problem for a conservation law, 
J. Math. Anal. Appl. 38 (1972), 33--41.

\item{[\refDafermostwo]} Dafermos C.M., 
Generalized characteristics in hyperbolic conservation laws: 
a study of the structure and the asymptotic behavior of solutions, 
in ``Nonlinear Analysis and Mechanics: Heriot-Watt symposium'', 
ed. R.J. Knops, Pitman, London, Vol. 1 (1977),  1--58.

\item{[\refDafermosthree]} Dafermos C., 
{\it Hyperbolic Conservation Laws in Continuum Physics,\/} 
Grundlehren Math. Wissen., Vol. 325, Springer Verlag, 2000. 

\item{[\refDLM]} Dal Maso G., LeFloch P.G., and Murat F., 
Definition and weak stability of nonconservative products, 
J. Math. Pures Appl. 74 (1995), 483--548. 

\item{[\refEG]} Evans L.C. and Gariepy R.F., 
{\it Measure Theory and Fine Properties of Functions,\/}  
Studies in Advanced Mathematics, CRC Press, 1992.  

\item{[\refFilippov]}  Filippov A.F., 
Differential equations with discontinuous right hand-side, 
Math USSR-Sb. 51 (1960), 99--128. 
English transl. in A.M.S. Transl., Ser. 2, 42, 199--231.

\item{[\refGoatinL]} P. Goatin and P.G. LeFloch, 
The sharp ${\bold L^1}$ continuous dependence
of solutions of bounded variation for hyperbolic 
systems of conservation laws, 
Arch. Rational Mech. Anal. (2001), to appear. 

\item{[\refHL]} Hu J.-X. and LeFloch P.G.,
$L^1$ continuous dependence for systems of conservation laws, 
Arch. Rational Mech. Anal. 151 (2000), 45--93. 

\item{[\refLax]} Lax P.D., 
Shock wave and entropy, in ``Contributions to Nonlinear 
Functional Analysis'', ed. E. Zarantonello, Acad. Press, New York, 1971, 
pp.~603--634. 

\item{[\refLeFlochone]} LeFloch P.G., 
An existence and uniqueness result for two nonstrictly hyperbolic systems, 
IMA Volumes in Math. and its Appl. 27,``Nonlinear 
evolution equations that change type'', ed.~B.L. Keyfitz and M. Shearer, 
Springer Verlag (1990), pp.~126--138. 

\item{[\refLeFlochtwo]} LeFloch P.G., 
An introduction to nonclassical shocks of systems of conservation laws, 
Proc. International School on Hyperbolic Problems, 
Freiburg, Germany, Oct. 97, D. Kr\"oner, 
M. Ohlberger and C. Rohde eds., Lect. Notes Comput. Eng., Vol.~5, 
Springer Verlag, 1998, pp.~28--72. 

\item{[\refLeFlochthree]} LeFloch P.G., 
Well-posedness theory for hyperbolic systems of conservation laws, 
to appear. 

\item{[\refLeFlochfour]} LeFloch P.G., 
{\it Hyperbolic Systems of Conservation Laws: 
The Theory of Classical and Nonclassical Shock Waves,\/}
Lecture notes, in preparation. 

\item{[\refLL]} LeFloch P.G. and Liu T.P.,
Existence theory for nonlinear hyperbolic systems in nonconservative form,
Forum Math. 5 (1993), 261--280.

\item{[\refLX]} LeFloch P.G. and Xin Z.P., 
Uniqueness via the adjoint problems for systems of conservation laws, 
Comm. Pure Appl. Math. 46 (1993), 1499--1533. 

\item{[\refLYone]} Liu T.P. and Yang T.,
A new entropy functional for scalar conservation laws, 
Comm. Pure Appl. Math. 52 (1999), 1427--1442. 

\item{[\refLYtwo]} Liu T.P. and Yang T.,
$L^1$ stability of weak solutions for 2x2 systems of 
hyperbolic conservation laws, 
J. Amer. Math. Soc. 12 (1999), 729--774. 

\item{[\refLYthree]} Liu T.P. and Yang T., 
Well-posedness theory for hyperbolic conservation laws,  
Comm. Pure Appl. Math. 52 (1999), 1553--1580.

\item{[\refVolpert]} Volpert A.I.,
The space BV and quasilinear equations,
Math. USSR Sbornik 73 (1967), 225--267.

\enddocument